\documentclass[generic,11pt,preprint]{imsart} 
\usepackage{amsmath,amsthm,amssymb,amsfonts} 
\RequirePackage[numbers]{natbib}
\RequirePackage[colorlinks,citecolor=blue,urlcolor=blue]{hyperref}
\allowdisplaybreaks

\startlocaldefs
\newtheorem{theorem}{Theorem}

\newtheorem{corollary}{Corollary}

\theoremstyle{definition}
\newtheorem{remark}{Remark}
\endlocaldefs

\begin{document}
\begin{frontmatter}

\title{Uncertainty quantification for robust variable selection and multiple testing}  
\runtitle{UQ for Variable selection and multiple testing}

\begin{aug}
\author{\fnms{Eduard} \snm{Belitser}$^1$\ead[label=e1]{e.n.belitser@vu.nl}}
\and
\author{\fnms{Nurzhan} \snm{Nurushev}$^2$\ead[label=e2]{nurushevn@gmail.com}}
\address{$^1$VU Amsterdam and\, $^2$Rabobank} 

\runauthor{E.\ Belitser and N.\ Nurushev}
\affiliation{$^1$VU Amsterdam and\;  $^2$Rabobank}
\end{aug}

\begin{abstract}
We study the problem of identifying the set of \emph{active} variables, termed in 
the literature as \emph{variable selection} or \emph{multiple hypothesis testing}, depending 
on the pursued criteria. For a general \emph{robust setting} of non-normal, possibly 
dependent observations and  a generalized notion of \emph{active set}, we propose a  
procedure that is used simultaneously for the both tasks, variable selection and multiple testing. 
The procedure is based on the \emph{risk hull minimization} method, but can also be obtained 
as a result of an empirical Bayes approach or a penalization strategy. 
We address its quality via various criteria: the Hamming risk, FDR, FPR, FWER, NDR, FNR, 
and various \emph{multiple testing risks}, e.g., MTR=FDR+NDR; 
and discuss a weak optimality of our results.
Finally, we introduce and study, for the first time, the \emph{uncertainty quantification problem} 
in the variable selection and multiple testing context in our robust setting.
\end{abstract}


\end{frontmatter}

\section{Introduction}

We are concerned with the problem of  identifying the set of \emph{active} 
(or \emph{significant}) variables. This task appears in a wide variety of applied 
fields as genomics, functional MRI, neuro-imaging, astrophysics, among others. 
Such data is typically available on a large number of observation units,
which may or may not contain a signal; the signal, when present, may be relatively  
faint and is dispersed across different observation units 
in an unknown fashion (i.e., the \emph{sparsity pattern} is unknown to the observer). 
 
A prototypical application is GWAS (\emph{genome-wide association studies}),  
where millions of genetic factors are examined for their potential influence 
on phenotypic traits. Although the number of tested genomic locations sometimes 
exceeds $10^5$ or even $10^6$, it is often believed that only a small set of 
genetic locations have tangible influences on the outcome of the disease 
or the trait of interest. This is well modeled by the stylized assumption of signal sparsity. 

Depending on pursued criteria the problem is termed in the literature 
as either \emph{variables selection}, (also termed as  \emph{recovery of the sparsity pattern}), 
or \emph{multiple testing problem}. Commonly, the problems of variable selection and multiple 
testing  are studied separately in the literature, although there are conceptual 
similarities and connections between them. In fact, a variable selection method determines 
the corresponding multiple testing procedure and vice versa, the difference lies merely 
in different criteria for inference procedures.

\subsection{The observations model and the context} 
Suppose we observe  a high-dimensional ($\mathbb{R}^n$-valued) 
vector $X=(X_1,\ldots,X_n) \sim \mathbb{P}_\theta$ such that 
$(X-\theta)/\sigma$ satisfies Condition \eqref{cond_nonnormal} (for some $\sigma>0$), 
where $\theta=(\theta_1,\ldots, \theta_n) \in \mathbb{R}^n$ is an unknown 
high-dimensional \emph{signal}.  
Actually, $ \mathbb{P}_\theta= \mathbb{P}^{n}_{\theta,\sigma}$, 
but we will omit the dependence on $n$ and $\sigma$ in the sequel.
In other words, we observe
\begin{align}
\label{model}
X_i=\theta_i+ \sigma \xi_i, \quad 
\quad i \in [n]\triangleq\{1,\ldots,n\},
\end{align}
where $\xi=(\xi_1,\ldots, \xi_n)\triangleq (X-\theta)/\sigma$ is the ``noise'' vector and $\sigma>0$ 
is the known ``noise intensity''. We emphasize that we pursue a general \emph{distribution-free} 
setting: the distribution of $\xi$ is arbitrary, satisfying only Condition \eqref{cond_nonnormal} 
below. The purpose of introducing $\sigma$ is that certain extra information can be converted 
into a smaller noise intensity $\sigma$, a ``more informative'' model. For example, suppose we 
originally observed $X_{ij}$ with $\operatorname{E}_\theta X_{ij}=\theta_i$ and 
$\operatorname{Var}_\theta(X_{ij})=1$, 
such that $(X_{ij},j\in [m])$ are independent for each $i\in[n]$, for some $m=m_n \to \infty$ 
as $n\to\infty$. By taking $X_i=\frac{1}{m_n}\sum_{j=1}^{m_n} X_{ij}$, we obtain the model 
\eqref{model} with $\sigma^2=\frac{1}{m_n} \to 0$ as $n\to \infty$.
 We are interested in non-asymptotic results, which imply asymptotic ones if needed. 
Possible asymptotic regime is high-dimensional setup $n \to \infty$, the leading case 
in the literature for high dimensional models. Another  possible 
asymptotics is $\sigma\to 0$, accompanying $n\to \infty$ or on its own.

The general goal is (for now, loosely formulated) to select 
the \emph{active} (or \emph{significant}) coordinates $I_*(\theta)\subseteq[n]$ 
of the signal $\theta$, based on the data $X$. In the sequel, we will need to properly 
formalize the notion of \emph{active set} $I_*(\theta)$. In particular, in this paper we 
let $I_*(\theta)$ be not necessarily the support of $\theta$, $S(\theta)=\{i\in[n]: \theta_i \not= 0\}$.
The main motivation for this is that we may want to qualify some relatively small (but non-zero) 
coordinates of $\theta$ as ``inactive'', with the threshold depending on the number of 
such coordinates. On the other hand, if the non-zero coordinates $\theta_i$ are 
allowed to be arbitrarily close (relative to the noise intensity $\sigma$) to zero, 
then it is clearly impossible to recover the signal support. So, even when relaxing 
the notion of active set, there are still principal limitations as no method should 
be able to distinguish between $|\theta_i|\asymp \sigma$ and $\theta_i=0$. 
These limitations will be quantified by establishing an appropriate lower bound.  
To make the problem feasible, one needs either to impose some kind of \emph{strong signal
condition} on $\theta$ (typically done in the literature on variable selection), or  
somehow adjust (relax) the criterion that measures the procedure quality (typically done 
in the literature on the multiple testing). For example, certain procedures can control 
more tolerant criteria like FDR or NDR without any condition, but, as we show below, 
their sum can be controlled again only under some strong signal condition.

\subsection{Variable selection and multiple testing in the literature}

The standard, most studied situation in the literature is the particular case 
of \eqref{model}:
\begin{align}
\label{standard_model}
X_i \overset{\rm ind}{\sim} \mathrm{N}(\theta_i,\sigma^2), \;\;  i\in[n]; \quad
I_*(\theta) = S(\theta) 
\triangleq\{i\in[n]: \theta_i \not=0\},
\end{align} 
where the support $S(\theta)$ plays the role of active coordinates of $\theta$, and $\theta$ is 
assumed to be \emph{sparse} in the sense that $\theta\in \ell_0[s]=\{\theta\in\mathbb{R}^n: 
|S(\theta)| \le s\}$ with $s=s_n=o(n)$ as $n\to \infty$.

Considering the situation \eqref{standard_model} for now,  
there is a huge literature on the active set recovery problem 
studied from various perspectives. Let $\hat{I}= \hat{I}(X)\subseteq[n]$ be a data 
dependent selector of the active set, $\eta_I= (\eta_i(I),i\in[n])=\big(1\{i\in I\}, i\in[n]\big)$ 
be the binary representation of  $I\subseteq [n]$. The historically first approach is via 
the probability of wrong recovery $\mathbb{P}_\theta(\hat{I} \not = S(\theta))$. 
The Hamming distance between $I_1$ and $I_2$ is $|\eta_{I_1} -\eta_{I_2}|\triangleq\sum_{i\in[n]}
|\eta_i(I_1)-\eta_i(I_2)| =|I_1\backslash I_2|+|I_2\backslash I_1|$.
One common measure of the quality of $\hat{I}$ is the expected Hamming loss 
(which we will call \emph{Hamming risk})  
\vspace{2mm}
\[
R_H(\hat{I},I_*)=\operatorname{E}_\theta |\eta_{\hat{I}}- \eta_{I_*}|
=\operatorname{E}_\theta  |\hat{I} \backslash I_*|+ \mathrm{E}_\theta |I_*\backslash \hat{I}|
=R_{FP}(\hat{I},I_*)+R_{FN}(\hat{I},I_*),
\] 
where $R_{FP}$ and $R_{FN}$ are the \emph{false positives} and \emph{false negatives} 
terms (in a way, Type I and Type II errors), respectively.
Note that $\mathrm{P}_\theta \big(\hat{I} \not = I_*\big)=
\mathrm{P}_\theta\big(|\eta_{\hat{I}}- \eta_{I_*}| \ge 1\big)\le
\mathrm{E}_\theta |\eta_{\hat{I}}-\eta_{I_*}|= R_H(\hat{I},I_*)$,
which means that the approach based on the Hamming risk provides stronger results. 

As is already well understood in many papers in related situations, in order to be able 
to recover $S(\theta)$, the non-zero signal $\theta_{S(\theta)}=(\theta_i, i\in S(\theta))$ 
has to satisfy some sort of \emph{strong signal condition}. If $\theta \in \ell_0[s]$ with 
polynomial sparsity parametrization $s = n^\beta$, $\beta\in (0,1)$ in the normal model 
\eqref{standard_model}, \cite{Donoho&Jin:2004} and \cite{Arias-Castro&Chen:2017} 
(see further references therein) express this condition in the form of the right scaling 
for the active coordinates $\theta_i^2 \gtrsim  \sigma^2 \log(n)$, $i\in S(\theta)$,
for the signal to be detectable. Actually, \cite{Donoho&Jin:2004} studied an idealized 
chi-squared model $Y_i\overset{\rm ind}{\sim} \chi^2_\nu(\lambda_i)$, $i\in[n]$, 
where $ \chi^2_\nu(\lambda_i)$ is a chi-square distributed random 
variable with $\nu$ degrees of freedom and non-centrality parameter $\lambda_i$. 
Squaring the both sides of \eqref{standard_model}, we arrive at 
the above chi-squared model with non-centrality parameters $\lambda_i = \theta_i^2$ and 
degree-of-freedom parameter $\nu=1$. 

Polynomial sparsity $s = n^\beta$ has been well investigated, especially in the normal 
model \eqref{standard_model}, with active coordinates 
$\theta_i^2 \gtrsim \sigma^2 \log(n)$, $i\in S(\theta)$. The situation with arbitrary 
signal sparsity $s$  is more complex:  
the right scaling for the active coordinates of $\theta\in\ell_0[s]$ becomes 
essentially (assuming $s=s_n \to \infty$ as $n \to \infty$) 
$\theta_i^2 \gtrsim \sigma^2 \log(n/s)$, $i\in S(\theta)$.
A  recent important reference on this topic is  \cite{Butucea&Ndaoud&Stepanova&Tsybakov:2018},
see also \cite{Ingster&Stepanova:2014}, \cite{Butucea&Stepanova:2017}. 
More on this is in Section \ref{sec_discussion}.

Inference on the active set can also be looked at from 
the multiple testing perspective. In classical multiple testing problem for the situation 
\eqref{standard_model}, one considers the following sequence of tests: 
\vspace{2mm}
\[
H_{0,i}: \theta_i=0 \quad  
\text{versus}   \quad H_{1,i}: \theta_i\neq 0,  \quad i\in[n]. 
\]
To connect to variable selection, notice that in multiple 
testing language, a variable selector $\hat{I}$ gives the multiple testing procedure which rejects 
the corresponding null hypothesis $H_{0,i}$, $i\in\hat{I}$, whereas $I_*(\theta)$  
encodes which null hypothesis do not hold. Typically, in multiple testing framework, 
one is interested in controlling (up to some prescribed level) of a type I error.
The most popular one is the so called \emph{False Discovery Rate} (FDR): 
$\text{FDR}(\hat{I}, I_*)=\operatorname{E}_\theta\frac{|\hat{I}\backslash I_*|}{|\hat{I}|}$ 
(with the convention $0/0=0$), the averaged proportion
of errors among the selected variables; in multiple testing terminology, 
the expected ratio of incorrect rejections to total rejections. 
This criterion, introduced in \cite{Benjamini&Hochberg:1995}, 
has become very popular because it is ``tolerant'' and ``scalable'' 
in the sense that the more rejections are possible, the more false positives are allowed.
It also delivers an adaptive signal estimator, see 
\cite{Abramovich&Benjamini&Donoho&Johnstone:2006}.

Besides FDR, we will study other 
known multiple testing criteria: \emph{Non-Discovery Rate} (NDR), 
\emph{k-Family-wise Error Rate} (k-FWER) introduced by \cite{Lehmann&Romano:2005}, 
\emph{False Non-Discovery Rate} (FNR) introduced 
by \cite{Genovese&Wasserman:2002} and \emph{False Positive Rate} (FPR). 
The FDR,  k-FWER and FPR have the flavor of type I error as they deal with the control of false positives, 
NDP and FNR have the flavor of type II error as they deal with the control of false negatives.
A weak optimality of our results in relation to \cite{Butucea&Ndaoud&Stepanova&Tsybakov:2018} 
is discussed in Section \ref{sec_discussion}.

In the multiple testing setting, most of theoretical studies rely on the fact that the null distribution 
is exactly known. In practice, it is often unreasonable to assume this, instead 
the null distribution is commonly (e.g., in genomics) implicitly defined as the ”background noise” 
of the measurements, and is adjusted via some pre-processing steps.
The issue of finding an appropriate null distribution has been popularized by a series of papers 
by Efron, (see \cite{Efron:2008}-\cite{Efron:2009} and further references therein), where the
concept of \emph{empirical null distribution} was introduced. 
A recent reference on this topic is \cite{Roquaine&Verzelen:2020}, see also further 
references therein. Our robust setting actually aligns well with the fact that the null hypothesis 
distribution is unknown, in fact, we avoid the problem of estimating the null distribution and 
obtain results that are robust over a certain (rather general) family of null distributions.

\subsection{Multiple testing risk, strong signal condition again}
 
In multiple testing setting, controlling type I error only is clearly not enough 
to characterize the quality of the procedures. For example, taking 
$\hat{I}=\varnothing$ gives the perfect FDR control: 
$\text{FDR}(\varnothing,I_*)=0$, but this is of course an unreasonable procedure. 
One needs  to control also some type II error, for example, the so called 
\emph{Non-Discovery Rate} (NDR) $\text{NDR}(\hat{I},I_*)=
\operatorname{E}_\theta\frac{|I_*\backslash \hat{I}|}{|I_*|}$ 
(or, the \emph{false Non-Discovery Rate} (FNR)).
Again, if considered as the only criterion, the NDR can be easily controlled
simply by taking another unreasonable procedure $\hat{I}=[n]$, yielding 
$\text{NDR}([n],I_*)=0$.
Thus, it is relevant to control these errors together, 
see \cite{Arias-Castro&Chen:2017}, \cite{Rabinovich&Ramdas&Jordan&Wainwright:2017}, 
\cite{Salomond:2017}, \cite{Castillo&Roquain:2018} and further references in these papers.
For example, a criterion to control is  
the \emph{multiple testing risk} $\text{MTR}(\hat{I},I_*)=  
\text{FDR}(\hat{I},I_*)+ \text{NDR}(\hat{I},I_*)$.   
 Other choices of MTR are possible, as long as it is a combination of 
type I and type II errors. Apart from FDR+NDR, these are FDR+FNR, FPR+NDR and 
FPR+FNR. 

Relating the MTR to the Hamming risk $R_H$, notice that, the both MTR and $R_H$, 
although being different, always combine some sort of Types I and II errors, 
in essence controlling  the \emph{false positives} and \emph{false negatives} simultaneously. 
Clearly, the MTR is a more mild criterion as it is always a proportion: $\text{MTR} \ll R_H$.  

It is desirable that $\text{MTR}$ is as small as possible, 
e.g., converging to zero as $n\to \infty$. 
However, as is shown in related settings in \cite{Ingster&Stepanova:2014}, 
\cite{Arias-Castro&Chen:2017}, \cite{Rabinovich&Ramdas&Jordan&Wainwright:2017}, 
\cite{Salomond:2017}, \cite{Butucea&Ndaoud&Stepanova&Tsybakov:2018}, 
\cite{Castillo&Roquain:2018}, this is in general impossible, which is not surprising because  
the same kind of principal limitation occurs for the Hamming risk in the case of variable selection.  
Again some sort of \emph{strong signal condition} is unavoidable. 
The Hamming risk $R_H$ is a more severe quality measure than MTR, so
its convergence to zero should occur under a more severe strong/sparse 
signal condition. In some papers in related settings this difference is referred to as 
\emph{exact} and \emph{approximate} recovery of the active set.  
Yet another type of recovery, the so called \emph{almost full recovery} is studied in \cite{Butucea&Ndaoud&Stepanova&Tsybakov:2018}, this is 
the convergence of $R_H/|I_*|$ to zero. Below we introduce all these notions 
more precisely.

\subsection{The scope}
\label{sec_scope}
In this paper, we generalize the standard normal setting \eqref{standard_model} to 
a more general setting \eqref{model} in that we pursue the \emph{robust} 
inference in the sense that the distribution of the error vector $\xi$ is unknown, but assumed 
to satisfy only certain condition, Condition \eqref{cond_nonnormal} below. 
In particular, the $\xi_i$'s can be non-normal, not identically distributed, of non-zero mean, 
and even dependent; their distribution may depend on $\theta$.   

Next, we generalize the notion of active set $I_*(\theta)$ that is now not necessarily 
the sparsity class $\ell_0(s)$ and not necessarily with all large non-zero coordinates.
We propose to express strong signal conditions as a scale of classes $\{\Theta(K), K>0\}$   
for the signal $\theta$, for the both problems simultaneously, 
variable selection and multiple testing problems. 
For the problem of determining the active set to be well defined, the parameter 
$\theta$ has to possess a \emph{distinct} set $I_*(\theta)$ of active coordinates,
which is ensured by the condition $\theta\in\Theta(K)$ for some (sufficiently large) 
$K>0$. The sparsity is expressed by $|I_*(\theta)|$ and the strong signal 
condition by the fact that $\theta \in \Theta(K)$ for sufficiently large $K$ (depending on the 
constants from Condition \eqref{cond_nonnormal}). Varying the ``goodness'' of $\theta$ 
(combination \emph{sparsity}/\emph{strong signal}) exposes the so called 
\emph{phase transition} effect, separating the impossibility and possibility to recover the active set,
shortly discussed in Section \ref{subsec_ phase_transition}.  
 
Finally, in this paper we address the new problem of \emph{uncertainty quantification} (UQ) 
for the active set $\eta_*=\eta_*(\theta)=\eta(I_*(\theta))$, this is to be distinguished 
from the uncertainty quantification for the parameter $\theta$.
For the Hamming loss $|\cdot|$ on $\{0,1\}^n$, a confidence ball is 
$B(\hat{\eta},\hat{r})=\{\eta \in \{0,1\}^n: |\hat{\eta}-\eta| \le \hat{r}\}$,
where the center $\hat{\eta}=\hat{\eta}(X):\mathbb{R}^n \mapsto \{0,1\}^n$ and radius 
$\hat{r}=\hat{r}(X): \mathbb{R}^n \mapsto \mathbb{R}_+ =[0,+\infty]$ are
measurable functions of the data $X$. The goal is to construct such a confidence 
ball $B(\hat{\eta},C\hat{r})$ that for any $\alpha_1,\alpha_2\in(0,1]$ and some function 
$r(\eta_*(\theta))$, $r:\mathbb{R}^n\rightarrow \mathbb{R}_+$, there exist $C, c > 0$ such that
\begin{align}
\label{defconfball}
\sup_{\theta\in\Theta_{\rm cov}}\mathbb{P}_\theta\big(\eta_*(\theta)\notin 
B(\hat{\eta},C\hat{r})\big)\le\alpha_1, \quad
\sup_{\theta\in\Theta_{\rm size}}\mathbb{P}_\theta\big(\hat{r}\ge 
c r(\eta_*(\theta))\big)\le\alpha_2,
\end{align}
for some $\Theta_{\rm cov}, \Theta_{\rm size} \subseteq \mathbb{R}^n$. 
The function $r(\eta_*(\theta))$, called the \emph{radial rate}, is a benchmark 
for the effective radius of the confidence ball $B(\hat{\eta},C\hat{r})$. 
The first expression in \eqref{defconfball} is called \emph{coverage relation} 
and the second \emph{size relation}. To the best of our knowledge, 
there are no results on uncertainty quantification with the Hamming loss 
\eqref{defconfball} for the active set $\eta_*(\theta)$. 
It is desirable to find the smallest  $r(\eta_*(\theta))$, 
the biggest  $\Theta_{\rm cov}$ and $\Theta_{\rm size}$ such 
that \eqref{defconfball} holds and $r(\eta_*(\theta))\asymp 
R(\Theta_{\rm size})$, where $R(\Theta_{\rm size})$ is the optimal rate 
in estimation problem for $\eta_*(\theta)$. We derive some UQ results  
for the proposed selector $\hat{I}$. 

Typically, the so called \emph{deceptiveness} issue is pertinent to 
the UQ problem, meaning that the confidence set of the optimal size and 
high coverage can only be constructed for non-deceptive parameters 
(in particular, $\Theta_{\rm cov}$ cannot be the whole set $\mathbb{R}^n$). 
Being non-deceptive is expressed by imposing some condition on the parameter; 
for example, the EBR (excessive bias restriction) condition $\Theta_{\rm cov}=\Theta_{\rm ebr}\subset \mathbb{R}^n$, see \cite{Belitser:2017}--\cite{Belitser&Nurushev:2020}.
 Interestingly, there is no deceptiveness issue as such for our UQ problem. 
An intuition behind this is as follows: the problem of active set recovery is 
more difficult than the UQ-problem in a sense that solving the former problem implies 
solving the latter. Then the condition $\theta \in \Theta(K)$ for the parameter 
to have distinct active coordinates implies also that the parameter 
is non-deceptive. In our case, we will have $\Theta_{\rm cov}=
\Theta_{\rm size} = \Theta(K)$ for some $K >0$.

\subsection{Organization of the rest of the paper}
In Section \ref{section_preliminaries} we introduce some notation, the generalized notion 
of active coordinates, and describe the criteria and procedures for variable selection and 
multiple testing. In Section \ref{main_results} we  present the main results of the paper. 
In Section \ref{sec_discussion} we shortly discuss a weak optimality of our results 
and a phase transition effect. The proofs of the theorems are collected in 
Section \ref{proofs_theorems}, which, despite generality of the setting and results, we could keep 
completely self-contained and relatively compact.

\section{Preliminaries}
\label{section_preliminaries}

\subsection{Notation}
Denote the  probability measure of $X$ from the model (\ref{model}) by $\mathbb{P}_\theta$, 
and by $\operatorname{E}_\theta$ the corresponding expectation. 
For the notational simplicity, we skip the dependence on $\sigma$ and $n$ of these quantities 
and many others. Denote by $1\{s\in S\}=1_S(s)$ the indicator function of the set $S$, 
by $|S|$ the cardinality of the set $S$, the difference of sets 
$S\backslash S_0 =\{s \in S:\,  s\not\in S_0\}$. Let $[k]=\{1,\ldots,k\}$ and $[k]_0=\{0\} \cup [k]$ 
for $k\in\mathbb{N}=\{1,2,\ldots\}$. For $I\subseteq [n]$ define $I^c=[n] \backslash I$. 
If random quantities  appear in a relation, this relation should be understood 
in $\mathbb{P}_{\theta}$-almost sure sense.   
For two nonnegative sequences $(a_l)$ and $(b_l)$, $a_l \lesssim b_l$ means $a_l \le cb_l$ 
for all $l$ (its range should be clear from the context)
with some absolute $c>0$, and $a_l \asymp b_l$ means that $a_l \lesssim b_l$ and 
$b_l \lesssim a_l$. The symbol $\triangleq$ will refer to equality by definition, 
$\Phi(x) = \mathbb{P}(Z \le x)$ for $Z \sim \mathrm{N}(0,1)$.
Throughout we assume the conventions: $|\varnothing|=0$, $\sum_{I\in\varnothing}a_I=0$ 
for any $a_I\in\mathbb{R}$ and $0\log(a/0)=0$ (hence $(a/0)^0=1$) for any $a>0$, 
in all the definitions whenever $0/0$ occurs we set by default $0/0 =0$.
Introduce the function $\ell(x)=\ell_q(x)\triangleq x\log(qn/x)$, $x,q> 0$, increasing 
in $x\in [0,n]$ for all $q\ge e$.
Finally, introduce the notation of ordered $\theta_{1}^2,\ldots, \theta_{n}^2$:
$\theta_{[1]}^2\ge\theta_{[2]}^2\ge\ldots \ge \theta_{[n]}^2$,  
define additionally  $\theta_{[0]}^2=\infty$ and $\theta_{[n+1]}^2=0$.

Recall the binary representation of an $I\subseteq[n]$:
\[
\eta_I=(\eta_i(I), i\in[n])=(\mathrm{1}\{i\in I\}, i\in[n]) \in \{0,1\}^n.
\] 
Let $\check{\eta}=\eta_{\check{I}}$ be some data dependent selector 
(for some measurable $\check{I}=\check{I}(X)$) which is supposed to estimate
\[
\eta_*=\eta_*(\theta)=\eta_{I_*}=(1\{i\in I_*\}, i\in[n]), 
\]
for some ``true''  \emph{active set} $I_*=I_*(\theta)$. 
The Hamming distance between $\check{\eta}$ and $\eta_*$ is determined by \
the number of positions at which $\check{\eta}$ and $\eta_*$ differ:
\begin{align}
\label{def_hamming}
 |\check{\eta}-\eta_*|= \sum_{i=1}^n |\check{\eta}_i-\eta_i(I_*)|
=\sum_{i=1}^n \mathrm{1}\{\check{\eta}_i\neq\eta_i(I_*)\}
=|\check{I}\backslash I_*|+|I_*\backslash\check{I}|.
\end{align}

Now we give some definitions for the multiple testing framework. For a variable selector 
$\check{I}=\check{I}(X)\subseteq [n]$ (which is seen now as multiple testing procedure)
and the active set $I_*$ (which is seen now as the set of the true null hypothesis), 
introduce the quantities characterizing the quality of  the multiple testing procedure $\check{I}$.
The convention $0/0=0$ is used in the following definitions. 
The  \emph{false discovery proportion} (FDP) and \emph{false discovery rate} (FDR) are
\begin{align*}
\text{FDP}(\check{I})=\text{FDP}(\check{I},I_*) 
=\frac{|\check{I}\backslash I_*|}{|\check{I}|},\quad  
\text{FDR}(\check{I})=\text{FDR}(\check{I},I_*)=\operatorname{E}_\theta\text{FDP}(\check{I}).
\end{align*}
The  \emph{false positive proportion} and \emph{false positive rate} are
\begin{align*}
\text{FPP}(\check{I})=\text{FPP}(\check{I},I_*) 
=\frac{|\check{I}\backslash I_*|}{n-|I_*|},\quad  
\text{FPR}(\check{I})=\text{FPR}(\check{I},I_*)=\operatorname{E}_\theta\text{FPP}(\check{I}).
\end{align*}
The \emph{non-discovery proportion}  and \emph{non-discovery rate} are
\begin{align*}
\text{NDP}(\check{I})=\text{NDP}(\check{I},I_*) 
=\frac{|I_*\backslash \check{I}|}{|I_*|},\quad 
\text{NDR}(\check{I})=\text{NDR}(\check{I},I_*)=\operatorname{E}_\theta\text{NDP}(\check{I}).
\end{align*}
The \emph{false non-discovery proportion}  and \emph{false non-discovery rate} are
\begin{align*}
\text{FNP}(\check{I})=\text{FNP}(\check{I},I_*)  
=\frac{|I_*\backslash \check{I}|}{n-|\check{I}|},\quad 
\text{FNR}(\check{I})=\text{FNR}(\check{I},I_*)=\operatorname{E}_\theta\text{FNP}(\check{I}).
\end{align*}

Introduce the \emph{multiple testing risks} (MTR)  
as all possible sums of the probabilities of type I and type II errors.
The first $\mathrm{MTR}$ is  the sum of the FDR and the NDR:
\begin{align*}
\text{MTR}_1(\check{I})=\text{MTR}_1(\check{I},I_*)=\text{FDR}(\check{I},I_*)
+\text{NDR}(\check{I},I_*).
\end{align*}
The other MTR's are defined similarly: they are always the sums of two rates (out of 4)
whose numerators must be the quantities $|\check{I}\backslash I_*|$ and $|I_*\backslash \check{I}|$.
Apart from $\text{MTR}_1$, these are $\text{MTR}_2
=\text{FDR}+\text{FNR}$, $\text{MTR}_3=\text{FPR}+\text{NDR}$ and 
$\text{MTR}_4=\text{FPR}+\text{FNR}$.

Finally introduce the  \emph{k-family-wise error} (\text{k-FWER}) and 
\emph{k-family-wise non-discovery} (\text{k-FWNR}) rates:
\begin{align*}
\text{k-FWER}(\check{I},I_*)=\mathbb{P}_\theta(|\check{I}\backslash I_*|\ge k),
\quad
\text{k-FWNR}(\check{I},I_*)=\mathbb{P}_\theta(|I_*\backslash \check{I}|\ge k).
\end{align*}
In multiple testing settings, the k-FWER is  the probability of rejecting at least $k$ true null hypotheses. 
The case $k=1$ reduces to the control of the usual FWER.

\subsection{Criterion for selecting active variables}

Consider for the moment the estimation problem of $\theta$ when we use the projection 
estimators $\hat{\theta}(I)= (X_i 1\{i\in I\}, i\in[n])$, $I \subseteq[n]$.  
The quadratic loss of $\hat{\theta}(I)$ gives a theoretical criterion  
\[
\mathcal{C}^{\rm th} _1(I)=\|\hat{\theta}(I) -\theta\|^2=\sum_{i\in I^c} \theta_i^2 
+ \sigma^2\sum_{i\in I}\xi_i^2.
\]
The best choice of $I$ would be the one minimizing $\mathcal{C}^{\rm th} _1(I)$. 
However, neither $\theta$ nor $\xi$ are observed. Substituting unbiased estimator 
$X_i^2 -\sigma^2$ instead of $\theta_i^2$, $i\in I^c$, leads to the quantity 
\[
\mathcal{C}^{\rm th}_2(I)=\sum_{i\in I^c} X_i^2 +\sigma^2 |I| + \sigma^2\sum_{i\in I}\xi_i^2
=\|X-\hat{\theta}(I)\|^2+\sigma^2 |I| + \sigma^2\sum_{i\in I}\xi_i^2
\]
to minimize with respect to $I\subseteq [n]$, however still not usable in view of the term
$\sigma^2\sum_{i\in I}\xi_i^2$.
If instead of $ \sigma^2\sum_{i\in I}\xi_i^2$, we use its expectation $\sigma^2 |I|$, we arrive 
essentially at Mallows's $C_p$-criterion (and AIC in the normal case) 
\[
\mathcal{C}_{\rm Mallows}(I)=\|X-\hat{\theta}(I)\|^2+2\sigma^2 |I|.
\] 
However, it is well known that 
the $C_p$-criterion leads to overfitting. An intuitive explanation is
that using the expectation $\sigma^2 |I|$ as penalty in $\mathcal{C}_{\rm Mallows}(I)$ 
is too optimistic to control oscillations of its stochastic counterpart $\sigma^2\sum_{i\in I}\xi_i^2$.

The next idea is to use some quantity $p(I)$ (instead of $|I|$) that majorizes 
$\sum_{i\in I}\xi_i^2$ in the more strict sense that for some $K, H_0,\alpha>0$
and all $M\ge 0$
\begin{align}
\label{main_ineq}
\mathbb{P}_\theta\Big(\sup_{I\subseteq[n]}\big(\sum_{i\in I}\xi_i^2- Kp(I)\big) \ge M\Big)  
\le H_0 e^{-\alpha M}.
\end{align}

\begin{remark} 
This idea is borrowed from the \emph{risk hull minimization method} developed by Golubev 
in several papers; see \cite{Cavalier&Golubev:2006} and references therein.
\end{remark}

Thus, in view \eqref{main_ineq}, using $K\sigma^2 p(I)$ instead of $\sigma^2\sum_{i\in I}\xi_i^2$,  
we obtain a more adequate criterion $\mathcal{C}_3(I)=\|X-\hat{\theta}(I)\|^2+
\sigma^2 |I|+K\sigma^2 p(I)$. Since typically $|I| \lesssim p(I)$, the second term can 
be absorbed into the third, and we finally derive the criterion 
\begin{align}
\label{criterion_C}
\mathcal{C}(I)=\|X-\hat{\theta}(I)\|^2+K\sigma^2 p(I),
\end{align}
for sufficiently large constant $K>0$. It remains to determine $p(I)$, 
preferably smallest possible, for which \eqref{main_ineq} holds. 
First we state an assumption.

{\sc Assumption (A1).}
For some $p_0(I)$ such that $p_0(I) \le C_\xi |I|$, for some $C_\xi>0$, $H_\xi,\alpha_\xi>0$ 
and any $M\ge 0$,  
\begin{align}
\label{cond_nonnormal}
\tag{A1}
\sup_{\theta\in\mathbb{R}^n}\mathbb{P}_\theta\Big(\sum_{i\in I} \xi_i^2 \ge p_0(I) +M\Big) 
\le H_\xi e^{-\alpha_\xi M}, \quad I\subseteq [n].
\end{align}
If the distribution of $\xi$ does not depend on $\theta$ 
(in some important specific cases), 
there is no supremum over $\theta\in\mathbb{R}^n$. 
For independent $\xi_i$'s, \eqref{cond_nonnormal} holds with $p_0(I) \asymp |I|$, 
so that  indeed $p_0(I) \le C_\xi |I|$, $I\subseteq[n]$, for some $C_\xi>0$. 
 
\begin{remark}
\label{rem_only_cond}
For fixed constants $C_\xi, H_\xi, \alpha_\xi$, one can think of \eqref{cond_nonnormal} with 
$p_0(I)=C_\xi |I|$ as description of a class of possible measures, 
in the sequel denoted by $\mathcal{P}_A$, see Remark \ref{rem_uniformity} below.
\end{remark}

\begin{remark}
Condition \eqref{cond_nonnormal} is of course satisfied for independent normals 
$\xi_i \overset{\rm ind}{\sim} \mathrm{N}(0,1)$ and for bounded (arbitrarily dependent) $\xi_i$'s.
Recall that the $\xi_i$'s are not necessarily of zero mean, but for normals 
Condition \eqref{cond_nonnormal} is the weakest if $\mathrm{E}_\theta \xi_i =0$.
In a way, Condition \eqref{cond_nonnormal} prevents too
much dependence, but it still allows some interesting cases of dependent
$\xi _{i}$'s. For example, one can show (in the same way as in \cite{Belitser&Nurushev:2020})
that this conditions is fulfilled for $\xi _{i}$'s that follow an autoregressive model
AR(1) with normal white noise.
\end{remark}

Let $\eta= \sup_{I\subseteq[n]}\big(\sum_{i\in I}\xi_i^2- C_\xi |I|
- \alpha_\xi^{-1}\big[|I|+ \log \tbinom{n}{|I|}\big]\big)$. 
By \eqref{cond_nonnormal} and the union bound, it is not difficult to derive 
\begin{align}
\mathbb{P}_\theta(\eta \ge M) &\le 
\sum_{I\subseteq[n]} \mathbb{P}_\theta 
\Big(\sum_{i\in I} \xi_i^2 \ge C_\xi |I| +\alpha_\xi^{-1}\big[|I|+\log \tbinom{n}{|I|}\big]+ M\Big) \notag\\
\label{relation_eta}
&\le 
H_\xi e^{-\alpha_\xi M} \sum_{I\subseteq[n]} e^{-|I|} \tbinom{n}{|I|}^{-1}
= H_\xi e^{-\alpha_\xi M} \sum_{k=0}^n e^{-k} 
\le H_0 e^{-\alpha_\xi M}.
\end{align}
where $H_0=H_\xi/(1-e^{-1})$. As $\tbinom{n}{k} \le (\frac{en}{k})^k$, $k\in [n]$, we have
\[
C_\xi |I| +\alpha_\xi^{-1}\big[|I|+\log \tbinom{n}{|I|}\big] \le 
(C_\xi +\alpha_\xi^{-1}) |I|+ \alpha_\xi^{-1} |I| \log \tfrac{en}{|I|}
=\alpha_\xi^{-1} |I| \log \big( \tfrac{q_\xi n}{|I|}\big),
\]
where $q_\xi = e^{C_\xi \alpha_\xi+2}$.
The last two displays imply the following relation (that we will need later): 
for appropriate $M_\xi>0$ and $q=e^2$
\begin{align}
\label{ineq_for_later}
\sup_{\theta\in\mathbb{R}^n} \sum_{I\subseteq[n]} \mathbb{P}_\theta 
\Big(\sum_{i\in I} \xi_i^2 \ge M_\xi |I| \log\big(\tfrac{q n}{|I|}\big) + M\Big) 
&\le H_0 e^{-\alpha_\xi M}.
\end{align} 
Although we will use the relation \eqref{ineq_for_later} only for $q=e^2$,
it is also implied by \eqref{cond_nonnormal} for any $q>1$ with appropriately chosen 
$M_\xi=M_\xi(q)$, certainly if $M_\xi(q) \ge \alpha_\xi^{-1}\big(\tfrac{\log q_\xi}{\log q}+1\big)$. 

In view of \eqref{relation_eta} with \eqref{ineq_for_later}, we see that under 
\eqref{cond_nonnormal}, the criterion \eqref{main_ineq} is satisfied with 
$p(I)=\ell_q(|I|)=|I| \log\big(\tfrac{q n}{|I|}\big)$.
According to \eqref{criterion_C}, 
this motivates the definition of the so called  \emph{preselector}
\begin{align}
\label{I_MAP}
\tilde{I}=\tilde{I}(K)=\arg\min_{I\subseteq[n]}\Big\{\sum_{i \in I^c}X_i^2 +K \sigma^2 p(I)\Big\},
\end{align}
where $p(I)=p_q(I)\triangleq\ell_q(|I|)=|I| \log\big(\tfrac{q n}{|I|}\big)$, 
for some $K>0$ and $q=e^2$. If $\tilde{I}$ is  not unique, take, say,  the one with 
the biggest sum $\sum_{i\in\tilde{I}}(n-i)$. 
\begin{remark}
Actually, an interesting interplay between constants $K$ and $q$ is possible, making certain 
constant in the proofs sharper. But we fix the second constant $q=e^2$ in \eqref{I_MAP} 
for the sake of mathematical succinctness.  
\end{remark}

Notice that if  $X_i^2 > K \sigma^2\log \big(\tfrac{q n}{|\tilde{I}|+1}\big)$ (with $q=e^2$), 
then $i \in \tilde{I}$. On the other hand, if $i \in \tilde{I}$, then 
\begin{align*}
X_i^2 &\ge K\sigma^2\big[\log (\tfrac{q n}{|\tilde{I}|})-(|\tilde{I}|-1)
\log\big(1+\tfrac{1}{|\tilde{I}|-1}\big)\big]\\
&\ge
K\sigma^2\big[\log(\tfrac{q n}{|\tilde{I}|})-\tfrac{|\tilde{I}|-1}{|\tilde{I}|}]
\ge K\sigma^2\big[\log(\tfrac{q n}{|\tilde{I}|})-1] = K\sigma^2\log(\tfrac{en}{|\tilde{I}|}).
\end{align*}

Next, define the \emph{selector} $\hat{\eta}=(\hat{\eta}_i, i\in[n])$ 
(and respectively $\hat{I}$) of significant coordinates as
\begin{align}
\label{threshold}
\hat{\eta}_i=\hat{\eta}_i(K)=\mathrm{1}\big\{X_i^2\ge
K \sigma^2\log \big(\tfrac{q n}{|\tilde{I}|}\big)\big\}, \;\;
\hat{I}=\hat{I}(K)=\{i\in[n]: \hat{\eta}_i=1\}.
\end{align}
Notice that the selector is always a subset of the preselector:
$\hat{I}(K)\subseteq\tilde{I}(K)$. 
 
\begin{remark}
We have already mentioned that this procedure can be related to   
the \emph{risk hull minimization} (RHM) method developed by Golubev.
In principle, (almost) the same procedures can 
be derived as a result of the empirical Bayes approach with appropriately chosen prior, 
or as a result of the penalization strategy with appropriately chosen penalty, see 
\cite{Belitser&Nurushev:2020}. It is interesting that several methodologies deliver akin procedures. 
\end{remark}

\subsection{The notion of active set}
\label{subsec_conditions}
Suppose we consider an arbitrary $\theta$ and would still like somehow divide  all the entries 
of $\theta$ into the groups of \emph{active} and \emph{inactive} coordinates. 
Clearly, as active group, the traditional support set $S(\theta)=\{i\in[n]: \theta_i \not=0\}$ 
is not sensible for arbitrary $\theta$, because nonzero but smallish coordinates $\theta_i$ 
should possibly be assigned to the inactive group.

For an arbitrary $\theta \in \mathbb{R}^n$, define the {\it active set} 
$I_*(A,\theta)=I_*(A,\theta,\sigma^2)$ as follows: with $q=e^2$,
\begin{align}
\label{oracle}
I_*(A,\theta)=\arg\!\min_{I\subseteq[n]}r^2(I,\theta), \;\; \text{where} \;\;
r^2(I,\theta)\triangleq\sum_{i\in I^c} \theta_i^2+A\sigma^2 |I| \log(\tfrac{qn}{|I|}) 
\end{align}
and $I_*(A,\theta)$ is with the smallest sum $\sum_{i \in I_*(A,\theta)} i $ if the minimum 
is not unique. By the definition \eqref{oracle}, \[
r^2(\theta)\triangleq r^2(I_*(A,\theta),\theta)\le r^2(I,\theta) \quad 
\text{for any} \quad I\subseteq[n].
\]  
This implies that 
\begin{align}
\label{property_I}
\text{if} \quad
\theta_i^2 \ge A\sigma^2\log \big(\tfrac{q n}{|I_*|+1}\big), \quad
\text{then} \quad i\in I_*\subseteq [n],
\end{align}
where we denoted for brevity $I_*=I_*(A,\theta)$. 
Conversely, 
if $i \in I_*$, then 
\begin{align}
\theta_i^2 &\ge A\sigma^2\big[\log (\tfrac{qn}{|I_*|})-(|I_*|-1)\log\big(1+\tfrac{1}{|I_*|-1}\big)\big]
\notag\\
\label{property_Io}
&\ge
A\sigma^2\big[\log(\tfrac{qn}{|I_*|})-\tfrac{|I_*|-1}{|I_*|}]
\ge A\sigma^2\big[\log(\tfrac{qn}{|I_*|})-1] = A\sigma^2\log(\tfrac{en}{|I_*|}).
\end{align}
We will use the last property later on. Also notice that $I_*$ depends on 
the product $A\sigma^2$ rather than just on $A$. If we consider the asymptotic regime
$\sigma^2\to 0$, it is instructive to fix the product $A\sigma^2$ so that $A\to \infty$, which 
can be interpreted as if $\theta$ satisfies  more and more stringent strong signal condition. 

\begin{remark} 
To get an idea what $I_*(A,\theta)$ means, suppose in \eqref{oracle} 
we had $A\sigma^2|I|\log(qn)$ instead of  $A\sigma^2 |I| \log(\tfrac{qn}{|I|})$. 
Then active coordinates would have been all $i\in[n]$ corresponding to 
``large'' $\theta_i^2 \ge A\sigma^2 \log(qn)$. 
The definition \eqref{oracle} does kind of the same, but the requirement for being active becomes 
slightly more lenient if there are more ``large'' coordinates. 
The function $x\log(qn/x)$ is increasing (in fact, for all $q\ge e$) for $x\in[1,n]$ slightly slower 
than $x\log(qn)$, creating the effect of ``borrowing strength'' via the number of active coordinates: 
the more such coordinates, the less stringent the property of being active becomes.
\end{remark}

The family $\mathcal{I}=\mathcal{I}(\theta)=\{I^{\rm vsp}_k(\theta),\, k\in[n]_0\}$, 
with $I^{\rm vsp}_k(\theta) = \{i\in[n]: \theta^2_i\ge \theta_{[k]}^2\}$, is called 
\emph{variable selection path} (VSP). It consists of at most $n+1$ embedded sets:
\[
\varnothing=I^{\rm vsp}_0(\theta)\subseteq I^{\rm vsp}_1(\theta)\subseteq \ldots 
\subseteq I^{\rm vsp}_n(\theta)=[n].
\] 
Clearly, $I_*(A,\theta)=I^{\rm vsp}_{|I_*(A,\theta)|}(\theta)$, $\theta\in\mathbb{R}^n$. 
The function $g(A)=|I_*(A,\theta)|: \mathbb{R}_+ 
\mapsto \{0\}\cup\mathbb{N}$ is a non-increasing right continuous  
step function taking values $|S(\theta)|,\ldots, 0$, as $A$ increases from 0 to infinity.  
If some of $\theta_{[k]}$ coincide, the corresponding sets $I^{\rm vsp}_k(\theta)$ 
in the variable selection path $\mathcal{I}$ merge. Accounting for these merges, notice that 
the true support $S(\theta)$ is the last set in the variable selection path $\mathcal{I}(\theta)$. 

\begin{remark}
\label{rem_properties}
We state some further properties of the active set $I_*$ and the variable selection path $\mathcal{I}$.
\begin{itemize}
\item[(a)]
The family $\{I_*(A,\theta), A\ge 0\}$ reproduces the variable selection path $\mathcal{I}$  
\[
\{I_*(A,\theta), A\ge 0\}=\mathcal{I}(\theta).
\]

\item[(b)]
For any $0\le A_1\le A_2$ and any $\theta\in\mathbb{R}^n$ , we have 
\[
\varnothing\subseteq I_*(A_2,\theta)\subseteq I_*(A_1,\theta)\subseteq S(\theta)\subseteq [n].
\] 

\item[(c)]
If for some $I\subseteq[n]$, 
$\theta_i^2\ge  A\sigma^2 \log(qn/|I|)$ for all $i\in I$ and 
$\theta_i^2\le  A\sigma^2\log(q)$ for all $i\in I^c$, then $I_*(A,\theta)=I$.
In particular, if $\theta_i^2\ge A\sigma^2 \log(qn/|S(\theta)|)$ for all 
$i \in S(\theta)$ and some $A>0$, then $I_*(A',\theta)=S(\theta)$ for all $A'\le A$. 
\end{itemize} 
\end{remark}
 
\section{Main results}
\label{main_results}
In this section we present the main results.

\subsection{Control of the preselector $\tilde{I}$}
First, we establish the results on  over-size and under-size control of 
the preselector $\tilde{I}=\tilde{I}(K)$  defined by \eqref{I_MAP}.
Recall that $\ell(x)=\ell_q(x)=x\log(\tfrac{qn}{x})$, $x\ge 0$, $q=e^2$.

\begin{theorem}
\label{theorem_1}
Let $\tilde{I}=\tilde{I}(K)$ be defined by \eqref{I_MAP}, $I_*(A,\theta)$ be defined 
by \eqref{oracle}, let $H_0$ be from 
\eqref{ineq_for_later}. Then for any $A_0$ there exist (sufficiently large) $K_0$ and 
constants $M_0,\alpha_0>0$ (depending on $A_0$) such that  for any $M\ge 0$,
\begin{align}
\tag{i}
\label{th1_i}
\sup_{\theta\in\mathbb{R}^n}\mathbb{P}_\theta
\big(\ell(|\tilde{I}(K_0)|)\ge  M_0 \ell(|I_*(A_0,\theta)|) +M\big)
\le H_0 e^{-\alpha_0 M}.
\end{align}
In particular, this implies that there exists $M_1>0$ such that for any $M\ge 0$
\begin{align}
\tag{i'}
\label{th1_i'}
\sup_{\theta\in\mathbb{R}^n}
e^{\alpha_0\ell(|I_*(A_0,\theta)|)}
\mathbb{P}_{\theta} \big(|\tilde{I}(K_0)|\ge M_1 |I_*(A_0,\theta)|+M \big)\le 
H_0 e^{-\alpha_0M/2}.
\end{align}

For any $K_1>0$, $\delta\in[0,1)$, there exist $A_1,\alpha_1,\alpha'_1>0$ 
(depending on $\delta,K_1$) such that for any $M\ge 0$
\begin{align}
\tag{ii}
\label{th1_ii}
\sup_{\theta\in\mathbb{R}^n}
e^{\alpha'_1 \ell(|I_*(A_1,\theta)|)} \mathbb{P}_\theta\big(\ell(|\tilde{I}(K_1)|)\le 
\delta \ell(|I_*(A_1,\theta)|)-M\big)
\le H_0 e^{-\alpha_1 M}.
\end{align}
In particular, this implies that there exist $A_1,\alpha'_1>0$ such that 
\begin{align}
\tag{ii'}
\label{th1_ii'}
\sup_{\theta\in\mathbb{R}^n} e^{\alpha'_1 \ell(|I_*(A_1,\theta)|)}
\mathbb{P}_\theta\big(|\tilde{I}(K_1)|\le  \delta |I_*(A_1,\theta)|\big)
\le H_0.
\end{align}
\end{theorem}

\begin{remark}
\label{rem_rem8}
A couple of remarks are in order.
\begin{itemize}
\item[(a)]
We can obtain another formulation of the property \eqref{th1_i} 
(and respectively  \eqref{th1_i'}): for any sufficiently large $K_0$ (e.g., $K_0 > 2M_\xi$) 
there exist $A_0>0$ and constants $M_0, \alpha_0>0$ such that for any $M\ge 0$, \eqref{th1_i} holds. 
\item[(b)]
We can also derive another formulation of  the property \eqref{th1_ii} 
(and respectively  \eqref{th1_ii'}): for any $A_1>0$ there exist (sufficiently small) 
$K_1>0$ and $\delta\in[0,1]$ such that for any $M\ge 0$, \eqref{th1_ii} holds.
 
\item[(c)]
Similarly to  \eqref{th1_i'}, we could establish \eqref{th1_ii'} in the following form:
\[
\sup_{\theta\in\mathbb{R}^n}
e^{\alpha'_1\ell(|I_*(A_1,\theta)|)}
\mathbb{P}_\theta \big(|\tilde{I}(K_1)|\le  \delta |I_*(A_1)| -M\big)
\le H_0 e^{-\alpha_1 M}.
\]
\end{itemize}
\end{remark}

From now on we fix some sufficiently large $K>0$ (such that, according to Remark 
\ref{rem_rem8}, there exists $A_0$ for which 
\eqref{th1_i} and \eqref{th1_i'} are fulfilled) and compute the corresponding preliminary selector 
$\tilde{I}(K)$. For this $K$, the properties \eqref{th1_i} and \eqref{th1_ii} 
(or, \eqref{th1_i'} and  \eqref{th1_ii'}) of Theorem \ref{theorem_1}
provide {\emph{separately} over-size and under-size control of $\tilde{I}$, 
with some $A_0(K)$ and $A_1(K)$, respectively. 
By analyzing the proof, we see that always $A_0(K) \le A_1(K)$. Hence  
$I_*(A_1,\theta) \subseteq I_*(A_0,\theta)$ (as it should  be),  
forming a \emph{shell} $I_*(A_0,\theta) \backslash I_*(A_1,\theta)$ in the VSP.
\begin{remark}
If we fix some distribution of $\xi$ satisfying Condition \eqref{cond_nonnormal},
another way of defining $K$ is to think of it as the smallest constant $K_0$ such that 
there exists $A_0(K_0)$ (by (a) of Remark \ref{rem_rem8}) for which \eqref{th1_i}  
is fulfilled for some $H_0\le H'_0<\infty$ and $\alpha_0\ge \alpha'_0>0$.
\end{remark}

\subsection{Set of signals with distinct active coordinates}
In the light of Theorem \ref{theorem_1}, we can say informally that $\tilde{I}$ 
``lives in an inflated shell''  between $I_*(A_1(K),\theta)$ and $I_*(A_0(K),\theta)$, 
which can be thought of as \emph{indifference zone} for the selector $\tilde{I}=\tilde{I}(K)$.
In general, the sets $I_*(A_1(K),\theta)\subseteq I_*(A_0(K),\theta)$ can be far apart, 
and $\tilde{I}$ may have too much room to vary. This means that 
the corresponding signal $\theta$ does not have distinct active and inactive coordinates, 
active coordinates as such are not identifiable.

The values of the constants $A_0(K)$, $A_1(K)$ evaluated in the proof of 
the above theorem are of course far from being optimal as we use rather rough bounds 
in the course of our argument. However, the main message of Theorem \ref{theorem_1} is that 
constants exist such that \eqref{th1_i'} and \eqref{th1_ii'} are fulfilled.  
This motivates the following definition of the \emph{set of signals with distinctive active 
and inactive coordinates}.

{\sc Definition.} Fix some $K,M_1,\delta >0$, and define $A_0(K)$ to be the biggest constant 
for which \eqref{th1_i'} holds with some $M'_1 \le M_1$ and $A_1(K)$ be the smallest 
constant for which \eqref{th1_ii'} is fulfilled with some $\delta' \ge \delta$. 
Introduce the set 
\begin{align}
\label{def_Theta(K)}
\Theta(K)=\{\theta\in\mathbb{R}^n: I_*(A_1(K), \theta) =I_*(A_0(K), \theta)\}.
\end{align}
In what follows, denote for brevity $I_*=I_*(A_1(K),\theta)$ and $I^*=I_*(A_0(K),\theta)$. 

Remember that the quantities $I_*$ and $I^*$  depend on $\theta$. 
The constants $A_0(K) = A_0(K,M_1,\delta)$ and $A_1(K) = A_1(K,M_1,\delta)$ 
(not  depending on $\theta$) exist in view of Theorem \ref{theorem_1}. 
This is the main mission of Theorem \ref{theorem_1}.
In essence, these constants are defined to be those which make the room between 
$I_*$ and $I^*$ (where $I_*=I_*(A_1(K),\theta) \subseteq I_*(A_0(K),\theta) =I^*$)
as small as possible uniformly over $\theta\in\mathbb{R}^n$. 
Hence, imposing $I_*=I^*$ determines a subset of $\mathbb{R}^n$ for which we can 
provide simultaneous control for oversizing and undersizing of $\tilde{I}$ by Theorem \ref{theorem_1}.
This means $\Theta(K)$ describes a set of signals with distinctive active and inactive 
coordinates: if $\theta\in\Theta(K)$, no indifference zone is allowed 
so that the active coordinates are well defined  as $I_*=I_*(A_1(K),\theta)$.

The above definition \eqref{def_Theta(K)} is still somewhat implicit, on the other hand 
it generalizes the traditional strong signal requirement. Indeed,  in view of property (c) 
from Remark \ref{rem_properties}, 
if $\theta_i^2 \ge A_1(K)\sigma^2 \log(\tfrac{qn}{|S(\theta)|})$ for all $i \in S(\theta)$, 
then $I_*=I_*(A_1(K),\theta)=S(\theta)=I_*(A_0(K),\theta)=I^*$ (in fact, $I_*(A,\theta)=S(\theta)$ 
for all $A\in [0,A_1(K)]$), so that $\theta\in\Theta(K)$.

\subsection{Results on active set recovery and multiple testing}
Now we establish the control of all the introduced quality measures 
(FPR, NDR, the Hamming rate, etc.) for the proposed 
active set selector $\hat{I}$. As consequence, we derive that  
the procedure  $\hat{I}$ matches (in a certain sense) the lower bound results 
from the previous section, establishing the optimality of the proposed procedure.
 
\begin{theorem}
\label{theorem_2}
Let $\hat{\eta}$ and $\hat{I}=\hat{I}(K)$ be  defined by \eqref{threshold}; 
$I^*=I_*(A_0(K))$, $I_*=I_*(A_1(K))$ and $\Theta(K)$ be defined by \eqref{def_Theta(K)} 
for sufficiently large $K>0$, and let $\eta_*=\eta_{I_*}$. Then there exist  constants $H_1,H_2,\alpha_2,\alpha_3>0$ such that, uniformly in $\theta\in\mathbb{R}^n$,  
\begin{align}
\label{th2_i}
\textnormal{FPR}(\hat{I},I^*)&=
\operatorname{E}_\theta \tfrac{|\hat{I}\backslash I^*|}{n-|I^*|}
\le
H_1 (\tfrac{n}{|I^*|\vee 1})^{-\alpha_2},\\
\label{th2_ii}
\textnormal{NDR}(\hat{I},I_*)&
=\operatorname{E}_\theta \tfrac{|I_*\backslash\hat{I}|}{|I_*|}
\le 
H_2 (\tfrac{n}{|I_*|})^{-\alpha_3}.
\end{align}
Moreover, there exist $H_3,\alpha_4>0$ such that, uniformly in $\theta\in \Theta(K)$,
\begin{align} 
\label{th2_iii} 
R_H(\hat{I},I_*)=\operatorname{E}_\theta | \hat{\eta} - \eta_*|=
\operatorname{E}_\theta \big(|\hat{I}\backslash I_*|+|I_*\backslash \hat{I}|\big)
&\le H_3 n (\tfrac{n}{|I_*|\vee 1})^{-\alpha_4}. 
\end{align} 
\end{theorem}

From the above theorem, the next corollary follows immediately.
It describes control of k-$\text{FWER}$, k-$\text{FWNR}$,
the probability of wrong discovery, and the so called \emph{almost full recovery} 
(relation \eqref{almost_full_recovery} below).
\begin{corollary}
\label{col1}
Uniformly in $\theta\in \mathbb{R}^n$,
\begin{align*}
\textnormal{k-FWER}(\hat{I}, I^*)
&=\mathbb{P}_\theta(|\hat{I}\backslash I^*|\ge k)
\le  
\tfrac{H_1}{k} (n-I^*)(\tfrac{n}{|I^*| \vee 1})^{-\alpha_2}, \\
\textnormal{k-FWNR}(\hat{I}, I_*)
&=\mathbb{P}_\theta(|I_*\backslash \hat{I}|\ge k)
\le  
\tfrac{H_2}{k} I_*(\tfrac{n}{|I_*|})^{-\alpha_3}. 
\end{align*}
Uniformly in $\theta\in \Theta(K)$, $\mathbb{P}_\theta(\hat{I}\neq I_*)\le 
H_3 n (\tfrac{n}{|I_*|})^{-\alpha_4}$ and 
\begin{align}
\label{almost_full_recovery}
\tfrac{R_H(\hat{I},I_*)}{|I_*|}=\tfrac{1}{|I_*|} \operatorname{E}_\theta \big(|\hat{I}\backslash I_*|
+|I_*\backslash \hat{I}|\big) \le H_3 (\tfrac{n}{|I_*|\vee 1})^{-(\alpha_4-1)}. 
\end{align}
\end{corollary}

The following result establishes the control of $\textnormal{FDR}(\hat{I})$ 
and $\textnormal{FNR}(\hat{I})$.
\begin{theorem}
\label{theorem_3}
With the same notation as in Theorem \ref{theorem_2}, 
there exist constants $H_5,H_6,\alpha_5,\alpha_6>0$ such that, 
uniformly in $\theta\in\Theta(K)$,  
\begin{align}
\label{th3_i}
\textnormal{FDR}(\hat{I},I_*)&=
\operatorname{E}_\theta \tfrac{|\hat{I}\backslash I_*|}{|\hat{I}|}
\le
H_5 (\tfrac{n}{|I_*|\vee 1})^{-\alpha_5},\\
\label{th3_ii}
\textnormal{FNR}(\hat{I},I_*)&=
\operatorname{E}_\theta \tfrac{|I_*\backslash \hat{I}|}{n-|\hat{I}|}
\le 
H_6 (\tfrac{n}{|I_*|})^{-\alpha_6}.
\end{align}
\end{theorem}

Theorems \ref{theorem_2} and \ref{theorem_3} imply the next corollary. 
\begin{corollary}
\label{col2}
For some constants $H_7,\alpha_7>0$, 
uniformly in $\theta\in \Theta(K)$,
\[
\mathrm{MTR}_l(\hat{I},I_*) \le 
H_7 (\tfrac{n}{|I_*| \vee 1})^{-\alpha_7}, \quad l=1,\ldots, 4.
\]
\end{corollary}

\begin{remark}
\label{rem_unif2}
In view of Remark \ref{rem_only_cond}, 
the results of Theorems \ref{theorem_1},  \ref{theorem_2} and 
Corollaries \ref{col1} and \ref{col2} hold also uniformly over 
all the measures $\mathbb{P}_\theta$ satisfying Condition 
\eqref{cond_nonnormal}.
\end{remark}

\begin{remark}
\label{rem_alpha}
Notice that all the powers $\alpha_i$'s in the above theorems and corollaries depend 
only on  $K$ (also via $A_0=A_0(K)$ and  $A_1=A_1(K)$) and the constants from Condition \eqref{cond_nonnormal}. 
Basically, the strong signal condition is reflected by the power $\alpha$'s: 
the stronger the signal, the bigger the $\alpha$. 
\end{remark}

Notice that $\textnormal{FPR}$ and \textnormal{NDP} are controlled uniformly in 
$\theta\in \mathbb{R}^n$, whereas all the other quantities only in
$\theta\in \Theta(K)$. As we already mentioned, the uniform control of 
either just Type I error or just Type II error is not much of a value, because 
this can always be achieved. It is a combination of the two types errors that 
one should try to control. The most natural choices of such combinations are 
the Hamming risk and the MTR's, studied in the present paper.
Another possible direction in obtaining interesting results is simultaneous control of 
Type I error (say, $\textnormal{FDR}$) and some estimation risk (or posterior 
convergence rate in case of Bayesian approach). Such a route is investigated in \cite{Castillo&Roquain:2018}. We should mention that $\hat{I}$ could also be derived 
as a result of empirical Bayes approach with appropriately chosen prior, 
and similar results could be derived on optimal estimation and posterior convergence rate.

Let us finally discuss possible asymptotic regimes. First, we note that asymptotics
$n\to \infty$ is not well defined, unless we describe how the true signal $\theta\in\mathbb{R}^n$ 
itself evolves with $n$. Assume that $\theta \in \mathbb{R}$
evolves with $n \in \mathbb{N}$ in such a way that  $\tfrac{|S(\theta)|}{n} \le p$ for some 
fixed $p\in [0,1)$.
Then from the definition \eqref{oracle} of active coordinates $I_*$, \eqref{property_I} and \eqref{property_Io}, it follows that $\frac{|I_*|}{n}\le p$, 
but it could also $\frac{|I_*|}{n}\to 0$.
This basically means that the signal is not getting ``less sparse'', in fact 
it can become ``more sparse'', making all the three criterions closer to zero. 
On the other hand, what can happen in this situation is that the signal 
is ``getting lost'' by spreading it over the bigger amount of coordinates. If, under 
growing dimension, we want the signal still to contain a certain portion of active 
coordinates, we need to make those coordinates more prominent, i.e., 
to strengthen the strong signal condition. The active coordinates should be 
increasing in magnitude when dimension is growing. 
This can also be attained by decreasing $\sigma^2$. Indeed, another observation is that 
if $\sigma^2 \to 0$, then $A\to \infty $ in the definition \eqref{oracle} of active coordinates $I_*$. 
One can interpret this as if the strong signal condition becomes more and more stringent. 
This in turn leads to $\alpha_4 \to \infty$ in \eqref{rel_optim1}.

\subsection{Quantifying uncertainty for the variable selector $\hat{\eta}$}  
Here we construct confidence ball $B(\hat{\eta},\hat{r})$ with optimal properties. 
Let $B(\hat{\eta},\hat{r})=\{\eta \in \{0,1\}^n: |\hat{\eta}-\eta| \le \hat{r}\}$, 
$\hat{\eta}=\hat{\eta}(K)$ and $\tilde{I}=\tilde{I}(K)$ be given by \eqref{I_MAP}. Define 
\begin{align}
\label{radius}
\hat{r}=\hat{r}(\tilde{I})= n (\tfrac{n}{|\tilde{I}|\vee 1})^{-\alpha'_4},
\end{align}
for some $\alpha'_4$ such that  $0<\alpha'_4< \alpha_4$, with $\alpha_4$  from 
Theorem \ref{theorem_2}.

The following theorem  describes 
the coverage and size properties of the confidence ball based on
$\hat{\eta}$  and $\hat{r}$. 
\begin{theorem}
\label{theorem_uncertainty}
With the same notation as in Theorem \ref{theorem_2}, let $\hat{r}$ be defined 
by \eqref{radius} and $ r_*=r_*(\theta)=n (\tfrac{n}{|I_*|\vee 1})^{-\alpha'_4}$. 
Then there exist constants $M'_1, H_7, H_8, \alpha_8, \alpha_{9}$ such that,
uniformly in $\theta\in\Theta(K)$,
\begin{align*}
\mathbb{P}_{\theta}\big(\eta_*\notin B(\hat{\eta},\hat{r})\big) 
&\le H_7 (\tfrac{n}{|I_*|\vee 1})^{-\alpha_8},\\
\mathbb{P}_{\theta}\big(\hat{r}\ge M'_1 r_* \big)
&\le H_8(\tfrac{n}{|I_*|\vee 1})^{-\alpha_{9}}.
\end{align*}
\end{theorem}
According to the UQ-framework \eqref{defconfball}, 
we have $\Theta_{\rm cov}=\Theta_{\rm size} = \Theta(K)$ for some $K >0$, 
$r(\eta_*(\theta)) =r_*(\theta)=n (\tfrac{n}{|I_*|\vee 1})^{-\alpha'_4}$. 
Notice that, according to \eqref{rel_optim1}, the radius $\hat{r}$ is \emph{optimal}, 
in a weak sense as it is up to the constant $\alpha'_4<\alpha_4$.

As we mentioned in Section \ref{sec_scope},  typically the so called \emph{deceptiveness} 
issue emerges in UQ problems. But in this case, interestingly, there is no deceptiveness 
issue as such for our UQ problem. A heuristic explanation is as follows: 
the problem of active set recovery is already  
more difficult than the UQ-problem in a sense that solving the former problem implies 
solving the latter. Basically, the condition $\theta \in \Theta(K)$ for the parameter 
to have distinct active coordinates implies also that the parameter 
is non-deceptive.

\section{Discussion: weak optimality, phase transition}
\label{sec_discussion}  

\subsection{Lower bounds} 
Define $\Theta_s(a)=\{ \theta\in \ell_0[s]:  |\theta_i| \ge a, i\in S(\theta)\}$.
Let $\Theta^+_s(a)$ be the version of $\Theta_s(a)$ 
when we put $\theta_i\ge a$ instead of $|\theta_i|\ge a$ in the definitions. Clearly,
$\Theta^+_s(a)\subset\Theta_s(a)$.  
For $\theta\in\Theta_s(a)$, the traditional active set is  
$I_*(\theta)=S(\theta)$. To ensure strict separation from the inactive set, 
one typically imposes $a\ge \bar{a}_n >0$ for appropriate $\bar{a}_n$. 

The minimax lower bound over the class $\Theta_s(a)$ for the problem of the recovery 
of the active set $I_*(\theta)=S(\theta)$ in the Hamming risk for the normal means 
model  was derived by \cite{Butucea&Ndaoud&Stepanova&Tsybakov:2018}.
Precisely, under the normality assumption $\xi_i \overset{\rm ind}{\sim} \mathrm{N}(0,1)$,
Theorem 2.2 from \cite{Butucea&Ndaoud&Stepanova&Tsybakov:2018}  states: 
for any $s<n$, $s'\in(0,s]$,
\[
r_H(\Theta^+_s(a)) \triangleq
\inf_{\check{\eta}} \sup_{\theta \in \Theta^+_s(a)}
\operatorname{E}_\theta |\check{\eta} -\eta_{S(\theta)}|
\ge s' \Psi_+(s,a) -4s' \exp\big\{-\tfrac{(s-s')^2}{2s}\big\},
\]
where 
$
\Psi_+(s,a) = (\tfrac{n}{s}-1) \Phi\big(-\tfrac{a}{2\sigma} - \tfrac{\sigma}{a}
\log(\tfrac{n}{s}-1) \big) +\Phi\big(-\tfrac{a}{2\sigma} +\tfrac{\sigma}{a}
\log(\tfrac{n}{s}-1) \big).
$
If $a^2 \le 2 \sigma^2 \log(\tfrac{n}{s} -1)$ then by taking $s'=s/2$ in the above display 
we get 
\begin{align}
\label{inconsistency}
r_H(\Theta^+_s(a)) \ge \tfrac{s}{2}\Phi(0) - 2se^{-s/8} =
s \big(\tfrac{1}{4} - 2e^{-s/8} \big)> 0.085 s
\end{align}
for $s\ge 20$. Expectedly, if $a^2 \le 2 \sigma^2 \log(\tfrac{n}{s} -1)$, 
it is impossible to achieve even consistency, so there is no point in considering this case. 
On the other hand, if $a^2> 2 \sigma^2 \log(\tfrac{n}{s} -1)$ and $n/s \ge 2.7$, then 
\[
\Phi\big(-\tfrac{a}{2\sigma} - \tfrac{\sigma}{a}
\log(\tfrac{n}{s}-1) \big) \ge 
\Phi(-a/\sigma) \ge (2/\pi)^{1/2} e^{-4 a^2/\sigma^2}.
\]
Assuming further  $\tfrac{a^2}{\sigma^2} \lesssim s$ 
(implying  $s\gtrsim \log n$) and taking again $s'=s/2$,  
\begin{align}
r_H(\Theta^+_s(a)) &\ge
\tfrac{(n-s)}{2}\Phi\big(-\tfrac{a}{2\sigma} - \tfrac{\sigma}{a}
\log(\tfrac{n}{s}-1) \big)- 2se^{-s/8}  \notag\\
\label{lower_bound1}
&\ge C_1(n-s)e^{-C_2 a^2/\sigma^2} - C_3 e^{-C_4 s} \ge 
C_5(n-s)e^{-C_6 a^2/\sigma^2}.
\end{align}
Assume that $a^2/\sigma^2 =A\log(\tfrac{en}{s})$ for some $A>2$
and $\log n \lesssim s \le n/2.7$. Then, under the normality assumption 
$\xi_i \overset{\rm ind}{\sim} \mathrm{N}(0,1)$, it follows from \eqref{lower_bound1} 
that for some $c_1, c_2>0$ (depending only on $A$)
\begin{align}
\inf_{\check{I}} \sup_{\theta \in \Theta_s(a)}
\operatorname{E}_\theta\big(|\check{I}\backslash S(\theta)|
&+|S(\theta)\backslash\check{I}|\big)
=\inf_{\check{I}} \sup_{\theta \in \Theta_s(a)} \operatorname{E}_\theta
|\eta_{\check{I}} - \eta_{S(\theta)}| \notag\\
\label{minimax_rate}
&=r_H(\Theta_s(a)) \ge r_H(\Theta^+_s(a))\ge c_1n (n/s)^{-c_2}.
\end{align}
In view of \eqref{inconsistency} and \eqref{minimax_rate}, even in the simplest 
normal model with $\theta \in \Theta_s(a)$, we need a strong signal condition $a \ge \bar{a}_n=A\sigma^2\log(\tfrac{en}{s})$ with $A>2$, just to avoid inconsistency 
in recovering the active set $S(\theta)$.
According to the terminology from \cite{Butucea&Ndaoud&Stepanova&Tsybakov:2018},
if $r_H(\Theta_s(a)) \to 0$ as $n\to\infty$, the \emph{exact recovery} of the 
active set $S(\theta)$ takes place; and \emph{almost full recovery} occurs if 
$r_H(\Theta_s(a))/s \to 0$ as $n\to\infty$ (assuming that $s>0$). 
The above lower bound \eqref{minimax_rate} reveals some sort of \emph{phase transition}. 
Indeed, the almost full recovery 
can occur if $s\ll n$ and the constant $A$ 
is sufficiently large, so that $c_2>1$, or, if $s\asymp n$ (but $s/n\le c<1$) and 
$c_2\to \infty$. The exact recovery is more 
difficult to fulfill, it can only occur if $n (n/s)^{-c_2}\to 0$ as $n\to\infty$. 
This is determined by the combination of two factors, the constant $c_2$ and 
the order of parameter $s$. The parameter $s$ describes the sparsity of the signal 
$\theta$ and the constant $c_2$ depends on $A$ which expresses the \emph{signal strength}.  

In view of property (c) from Remark \ref{rem_properties},
it follows that $\Theta_s(a(K)) \subseteq \Theta(K)$ for $a(K)=A_1(K)\sigma^2 \log(\tfrac{qn}{|S(\theta)|})$.
Property (b) from Remark \ref{rem_properties} implies also that 
$I_*= I_*(A_1(K), \theta)\subseteq S(\theta)$ for all $\theta\in\mathbb{R}^n$.
The last two facts and \eqref{minimax_rate} allow us to derive the following lower bound:  
\begin{align}
&\inf_{\check{I}} \sup_{\theta \in \Theta(K)}  
\operatorname{E}_\theta
\frac{|\check{I}\backslash I^*|+|I_*\backslash\check{I}|}{n (|I_*|/n)^{c_2}}
\ge 
\inf_{\check{I}} \sup_{\theta \in \Theta(K)} \operatorname{E}_\theta
\frac{|\check{I}\backslash I^*|+|I_*\backslash\check{I}|}{n(|S(\theta)|/n)^{c_2}} \notag \\
&\ge
\inf_{\check{I}} \sup_{\theta \in \Theta_s(a(K))} \operatorname{E}_\theta
\frac{|\check{I}\backslash S(\theta)|+|S(\theta)\backslash\check{I}|}{n(|S(\theta)|/n)^{c_2}}  
= 
\frac{r_H(\Theta_s(A_1(K)))}{n(s/n)^{c_2}}\ge c_1,
\label{lb_Theta}
\end{align}
where the distribution $\mathbb{P}_\theta$ is taken to be the product normal as 
in \eqref{standard_model}, $\log n \lesssim s \le n/2.7$ and $A_1(K)>2$. 
The normalizing factor for the Hamming risk is thus $n (|I_*|/n)^{c_2}$.

Recall that the bound \eqref{lb_Theta} is in the regime  
$A_1(K)>2$. Otherwise (i.e., when $A_1(K)\le2$), we have by \eqref{inconsistency} that 
\[
\inf_{\check{I}} \sup_{\theta \in \Theta(K)}
\operatorname{E}_\theta
\tfrac{|\check{I}\backslash I^*|+|I_*\backslash\check{I}|}{|I_*|}
\ge 0.085.
\] 
\begin{remark}
\label{rem_uniformity}
Notice that the measure $\mathbb{P}_\theta$ in the above lower bounds is
the product normal measure \eqref{standard_model}. However, 
the same lower bounds trivially hold when, instead of $\sup_{\theta\in\Theta(K)}$,
we take $\sup_{\mathbb{P}_\theta \in \mathcal{P}}$, where 
\[\mathcal{P}= \{ \mathbb{P}_\theta \in \mathcal{P}_A: \theta \in \Theta(K)\}, \quad 
\mathcal{P}_A=\{\mathbb{P}_\theta:  
\mathbb{P}_\theta \text{ satisfies \eqref{cond_nonnormal}}\}.
\] 
This is because this normal measure is one of many that satisfy \eqref{cond_nonnormal}.
Note that the constants $A_1(K), A_0(K)$ are then defined uniformly 
over  $\mathcal{P}_A$.
\end{remark}

\begin{remark}
Similarly to \eqref{lb_Theta}, it is easy to derive from \eqref{minimax_rate} that
$
\inf_{\check{I}} \sup_{\theta \in \Theta(K)}  
\frac{\mathrm{MTR}_l(\check{I},I_*)}{(|I_*|/n)^{c_2}} \ge c_1$, $l=1,\ldots, 4$.
But we conjecture that the right normalizing factor for the minimax 
MTR's should be $(|I_*|/n)^{c_2-1}$, this cannot be derived from \eqref{minimax_rate}.
\end{remark}


\subsection{Phase transition} 
\label{subsec_ phase_transition}  
Relating the lower bound \eqref{lb_Theta} with 
the results of the previous section, we claim that the selector $\hat{I}$ 
(and the corresponding $\hat{\eta}$) is \emph{optimal} in the following weak sense:
\begin{align}
\label{rel_optim1}
c_1 \le \inf_{\check{I}} \sup_{\theta \in \Theta(K)}   
\frac{R_H(\check{I},I_*)}{n (|I_*|/n)^{c_2}}, \;\; 
\sup_{\theta \in \Theta(K)}  
\frac{R_H(\hat{I},I_*)}{n ((|I_*| \vee 1)/n)^{\alpha_4}}
\le H_3.
\end{align}
In view of Remarks \ref{rem_uniformity}, the above relations hold also 
uniformly over all the measures $\mathbb{P}_\theta$ satisfying Condition 
\eqref{cond_nonnormal}.

Admittedly, the optimality in \eqref{rel_optim1} is very weak as it up to constants 
$c_2$, $\alpha_4$, which differ in general. But this is the best
we could achieve under the general robust setting of this paper. 
Constant $c_2$ in the lower bound is established for the normal 
submodel \eqref{standard_model} of our more general model \eqref{model}, whereas $\alpha_4$ 
is obtained uniformly for the general model \eqref{model} under Condition \eqref{cond_nonnormal} 
(thus determined by the constants from Condition \eqref{cond_nonnormal}). 

\begin{remark}
We should emphasize that if we want to match of upper and lower bounds we would need to  
specify the error distribution (or severely restrict the choice). This problem seemingly interesting 
and challenging does not align with the main focus of the present paper, the robust setting. 
\end{remark}

However, even these relatively loose lower and upper bounds in \eqref{rel_optim1}   
can demonstrate some sort of {\it phase transition} phenomenon, also in the general 
setting \eqref{model}. Precisely, the minimax Hamming $r_H$ risk (and MTR) can be either close to zero 
or not, depending on the combination of the signal sparsity $|I_*|$ and 
signal magnitude (how strong the signal is), reflected by the constants $c_2,\alpha_4$. 
The dependence of the normalizing factor on sparsity $|I_*|$ 
is only through the ratio $n/|I_*|$.
The ``informativeness'' of the model (how ``bad'' the noise is) plays a role as well in that it 
determines how large $K$ must be in the set $\Theta(K)$ for the upper bound in  \eqref{rel_optim1} 
to hold, which depends on the constants from Condition  \eqref{cond_nonnormal}.


\section{Proofs of  the theorems}
\label{proofs_theorems}
\begin{proof}[Proof of Theorem \ref{theorem_1}] 
For any $a,b\in\mathbb{R}$, $(a+b)^2\le 2a^2+2b^2$, 
hence also $-(a+b)^2 \le -a^2/2+b^2$. Using these elementary inequalities, 
the definition \eqref{I_MAP} of $\tilde{I}(K)$, 
we derive that, for any $I, I_0\subseteq[n]$,
\begin{align*}
&\mathbb{P}_{\theta}(\tilde{I}(K)= I) 
\le\mathbb{P}_{\theta}\Big(\sum_{i\in I^c} X_i^2+ K\sigma^2 \ell(|I|) \le 
\sum_{i\in I_0^c} X_i^2+K\sigma^2 \ell(|I_0|)  \Big) \notag\\
&= \mathbb{P}_{\theta}\Big(
\sum_{i\in I\backslash I_0} \tfrac{X_i^2}{\sigma^2} 
-\sum_{i\in I_0\backslash  I} \tfrac{X_i^2}{\sigma^2} \ge K[\ell(|I|)-\ell(|I_0|)] \Big)\notag\\
&\le \mathbb{P}_{\theta}\Big(
\sum_{i\in I\backslash I_0} (\tfrac{2\theta_i^2}{\sigma^2} +2\xi_i^2) 
-\sum_{i\in I_0\backslash  I} (\tfrac{\theta_i^2}{2\sigma^2}- \xi_i^2) 
\ge K(\ell(|I|)-\ell(|I_0|)) \Big) \notag\\
&= \mathbb{P}_{\theta}
\Big(\sum_{i\in I\backslash  I_0}2\xi_i^2 +\sum_{i\in I_0\backslash  I} \xi_i^2
\ge \sum_{i\in I_0\backslash  I}\tfrac{\theta_i^2}{2\sigma^2} 
-\sum_{i\in I\backslash I_0} \tfrac{2\theta_i^2}{\sigma^2} 
+K(\ell(|I|)-\ell(|I_0|)) \Big).
\end{align*}
In particular, for any $I, I_0$ such that $I_0 \subseteq I$, we have 
\begin{align}
\label{I-I_0}
\mathbb{P}_{\theta}(\tilde{I}(K)= I) \le 
\mathbb{P}_{\theta}\Big(
\sum_{i\in I\backslash  I_0}\xi_i^2 
\ge \tfrac{K}{2}(\ell(|I|)-\ell(|I_0|)) -\sum_{i\in I\backslash I_0} \tfrac{\theta_i^2}{\sigma^2}\Big),
\end{align}
and for any $I, I_0$ such that $I \subseteq I_0$, we have 
\begin{align}
\label{I_0-I}
\mathbb{P}_{\theta}(\tilde{I}(K)= I) \le 
\mathbb{P}_{\theta}\Big(\sum_{i\in I_0\backslash  I} \xi_i^2 
\ge \sum_{i\in I_0\backslash  I}\tfrac{\theta_i^2}{2\sigma^2} + K(\ell(|I|)-\ell(|I_0|))\Big).
\end{align}

Now we prove \eqref{th1_i}. 
For brevity, denote for now $I_*=I_*(A_0)=I_*(A_0,\theta)$.
If $A_0|I \backslash I_*|\log(\tfrac{qn}{|I\cup I_*|})<
\sum_{i\in I \backslash I_*}\tfrac{\theta_i^2}{\sigma^2}$ would hold for some 
$I\subseteq[n]$, then 
\begin{align*}
r^2_{A_0}(I\cup I_*,\theta) 
&= \sum_{i \not\in I\cup I_*} \theta^2_i+A_0\sigma^2 |I\cup I_*| 
\log(\tfrac{qn}{|I\cup I_*|})\\
&\le 
\sum_{i \not\in I\cup I_*} \theta^2_i +A_0\sigma^2 |I \backslash I_*|
\log(\tfrac{qn}{|I\cup I_*|})+A_0\sigma^2 |I_*|\log(\tfrac{qn}{|I_*|}) \\
&< 
\sum_{i \not\in I\cup I_*} \theta^2_i+\sum_{i\in I \backslash I_*} \theta_i^2
+ A_0\sigma^2 |I_*| \log(\tfrac{qn}{|I_*|}) \\
&=
\sum_{i \not\in I_*} \theta^2_i+A_0\sigma^2 |I_*| \log(\tfrac{qn}{|I_*|})=r^2_{A_0}(\theta),
\end{align*}
which contradicts the definition \eqref{oracle} of $I_*=I_*(A_0,\theta)$. Hence, 
\[
\sum_{i\in I \backslash I_*}\tfrac{\theta_i^2}{\sigma^2}\le 
A_0|I \backslash I_*|\log(\tfrac{qn}{|I\cup I_*|}) \le A_0 |I|\log(\tfrac{qn}{|I|})=A_0\ell(|I|).
\]
Using this and \eqref{I-I_0} with $I_0=I_*\cap I$ (so that $I\backslash  I_0=I\backslash  I_*$)
yields
\begin{align*}
\mathbb{P}_{\theta}(\tilde{I}(K_0)= I) &\le 
\mathbb{P}_{\theta}\Big(\sum_{i\in I\backslash  I_0}\xi_i^2 
\ge (\tfrac{K_0}{2}-A_0)\ell(|I|) - \tfrac{K_0}{2}\ell(|I_0|)\Big) \\
&\le \mathbb{P}_{\theta}\Big(\sum_{i\in I}\xi_i^2 
\ge (\tfrac{K_0}{2}-A_0)\ell(|I|) - \tfrac{K_0}{2}\ell(|I_*|)\Big).
\end{align*}

Let $\mathcal{J}=\{I\subseteq[n]: \ell(|I|) \ge M_0 \ell(|I_*(A_0)|)+M \}$, with $M_0=K_0/(2C_0)$ 
where $C_0=K_0/2-A_0-M_\xi>0$, which holds for any $K_0> 2(A_0+M_\xi)$. 
The last display and \eqref{ineq_for_later} imply that 
for any $\theta\in\mathbb{R}^n$
\begin{align*}
\mathbb{P}_{\theta} \big(\tilde{I}(K_0) \in \mathcal{J} \big)
&=
\sum_{I\in \mathcal{J}} \mathbb{P}_{\theta}(\tilde{I}(K_0)= I) \\&
\le 
\sum_{I\in \mathcal{J}} \mathbb{P}_{\theta}\Big(\sum_{i\in I}\xi_i^2 
\ge (\tfrac{K_0}{2}-A_0)\ell(|I|) - \tfrac{K_0}{2}\ell(|I_*|)\Big)\\
&= 
\sum_{I\in \mathcal{J}} \mathbb{P}_{\theta}\Big(\sum_{i\in I}\xi_i^2 
\ge M_\xi \ell(|I|)+C_0\ell(|I|) - \tfrac{K_0}{2}\ell(|I_*|)\Big)\\
&\le \sum_{I\subseteq[n]} \mathbb{P}_{\theta}\Big(\sum_{i\in I}\xi_i^2 
\ge  M_\xi \ell(|I|) + C_0M \Big) \le  H_0 e^{-\alpha_0 M},
\end{align*}
with $\alpha_0= C_0 \alpha_\xi$, 
thus ensuring the result of the assertion \eqref{th1_i} for any $\theta\in\mathbb{R}^n$.

To prove \eqref{th1_i'}, let $M'_1=4(M_0+1)$. If $|I|\ge M'_1|I_*|+M$ and $|I_*|\le qn/(M'_1)^2$, 
\begin{align*}
\ell(|I|) &=|I|\log\big(\tfrac{qn}{|I|}\big)\ge 
\tfrac{1}{2}\ell(|I|)+ \tfrac{1}{2}|I|\ge 
\tfrac{M'_1}{2}|I_*|\log\big(\tfrac{qn}{M'_1|I_*|}\big)+\tfrac{M}{2} \\
&\ge
\tfrac{M'_1}{4}|I_*| \log\big(\tfrac{qn}{|I_*|}\big) +\tfrac{M}{2}
=(M_0+1)\ell(|I_*|)+\tfrac{M}{2}.
\end{align*}
Hence, for $|I_*|\le qn/M_1^2$, by using  \eqref{th1_i},
\begin{align*}
\mathbb{P}_{\theta} \big(|\tilde{I}(K_0)|\ge M'_1 |I_*|+M\big) &\le 
\mathbb{P}_{\theta} \big(\ell(|\hat{I}(K_0)|)\ge (M_0+1)\ell(|I_*|)+\tfrac{M}{2}\big)\\
&\le
H_0 e^{-\alpha_0\ell(|I_*|) -\alpha_0M/2}.
\end{align*}
If $|I_*|> qn/(M'_1)^2$ and $M''_1=(M'_1)^2/q$, then we trivially obtain 
\[
\mathbb{P}_{\theta} \big(|\tilde{I}|\ge M''_1 |I_*|+M\big)
\le \mathbb{P}_{\theta} \big(|\tilde{I}|\ge M''_1 |I_*|\big)
\le \mathbb{P}_{\theta} \big(|\tilde{I}|>n \big)=0.
\]
Hence the choice $M_1=\max\{M'_1, M''_1\}$ ensures the second relation \eqref{th1_i'}.

Next, we prove the assertion \eqref{th1_ii}. For the rest of the proof, denote for brevity 
$I_*=I_*(A_1)$. Define $\mathcal{T}=\{I\in\mathcal{I}: \ell(|I|)\le \delta\ell(|I_*|)-M\}$, 
$\delta \in [0,1)$. If $I\in\mathcal{T}$, then $ \ell(|I|)\le \delta\ell(|I_*|) 
\le\ell(\delta |I_*|)$, implying that $|I| \le \delta |I_*|$, as $\ell(x)$ is increasing for $x\in[0,n]$. 
Hence, for any $I\in\mathcal{T}$,
\begin{align}
\label{relation_5}
\ell(|I_*\cup I|)&= |I_*\cup I| \log(\tfrac{qn}{|I_*\cup I|}) 
\le |I_*|\log(\tfrac{qn}{|I_*\cup I |})+ |I|\log(\tfrac{qn}{|I_*\cup I|})\notag\\
&\le (1+\delta) |I_*|\log(\tfrac{qn}{|I_*|}) -M= (1+\delta) \ell(|I_*|)-M.
\end{align}
Next, by \eqref{property_Io} and the fact that $|I| \le \delta |I_*|$, we obtain that 
for any $I\in\mathcal{T}$,
\begin{align}
\label{relation_4}
\sum_{i\in I_* \backslash I}\tfrac{\theta_{i}^2}{\sigma^2} &\ge 
|I_*\backslash I|A_1 \log(\tfrac{en}{|I_*|})
\ge A_1 (1-\delta)|I_*|\log(\tfrac{en}{|I_*|}) \ge \tfrac{A_1 (1-\delta)}{2}\ell(|I_*|)  ,
\end{align}

Denote for brevity $C_{A_1}=\tfrac{A_1(1-\delta)}{4}-K_1(1+\delta)$. Using the  relation 
\eqref{I_0-I} with $I_0=I_*\cup I$ (so that $I_0\backslash I=I_*\backslash I$), the relations 
\eqref{relation_5}, \eqref{relation_4}, and \eqref{ineq_for_later}, we derive that  
\begin{align}
\mathbb{P}_{\theta} & (\tilde{I}(K_1)\in\mathcal{T}) = \sum_{I \in \mathcal{T}}
\mathbb{P}_{\theta}(\tilde{I}(K_1)=I)
\notag \\
&\le \sum_{I \in\mathcal{T}} \mathbb{P}_{\theta}\Big(\sum_{i\in I_*\backslash  I} \xi_i^2 
\ge \sum_{i\in I_*\backslash  I}\tfrac{\theta_i^2}{2\sigma^2} + K_1(\ell(|I|)-\ell(|I_*\cup I|))\Big)
 \notag\\
&\le \sum_{I \in\mathcal{T}} \mathbb{P}_{\theta}\Big(\sum_{i\in I_*\backslash  I} \xi_i^2 
\ge \big( \tfrac{A_1(1-\delta)}{4} -K_1(1+\delta)\big)\ell(|I_*|)+K_1[\ell(|I|)+M]\Big) \notag\\
\label{rel6}
&= \sum_{I \in\mathcal{T}} \mathbb{P}_{\theta}\Big(\sum_{i\in I_*\backslash  I} \xi_i^2 
\ge  C_{A_1}\ell(|I_*|)+K_1\ell(|I|)+K_1M\Big)\\
&\le \sum_{I \in\mathcal{T}} \mathbb{P}_{\theta}\Big(\sum_{i\in I_*\backslash  I} \xi_i^2 
\ge M_\xi \ell(|I_*|) +(C_{A_1}-M_\xi)\ell(|I_*|)+K_1M \Big)   \notag\\
&\le \sum_{I\subseteq[n]} \mathbb{P}_{\theta}\Big(\sum_{i\in I} \xi_i^2 
\ge  M_\xi \ell(|I|) +(C_{A_1}-M_\xi)\ell(|I_*|) +K_1M\Big)  \notag\\
&\le H_0 e^{-\alpha_\xi M-\alpha_\xi(C_{A_1}-M_\xi) \ell(|I_*|)}
=H_0 e^{-\alpha_1 M-\alpha'_1\ell(|I_*|)},  \notag
\end{align}
where $\alpha_1=\alpha_\xi K_1$, $\alpha'_1= \alpha_\xi (C_{A_1}-M_\xi)$ and 
$A_1$ is assumed to be so large that $C_{A_1}>M_\xi$. This proves \eqref{th1_ii}.

Finally, we establish \eqref{th1_ii'}. Define 
$\mathcal{T}'=\{I\in\mathcal{I}: |I|\le \delta|I_*|\}$, $\delta \in [0,1)$.
The relations \eqref{relation_4}, \eqref{relation_5} are still valid for any $I\in\mathcal{T}'$.
As before, we obtain \eqref{rel6} (with $M=0$ and $\mathcal{T}'$ instead of $\mathcal{T}$), 
which we now continue as follows:
\begin{align*}
\mathbb{P}_{\theta}(\tilde{I}(K_1)\in\mathcal{T}')
&\le
\sum_{I \in\mathcal{T}'} \mathbb{P}_{\theta}\Big(\sum_{i\in I_*\backslash  I} \xi_i^2 
\ge  C_{A_1}\ell(|I_*|)+K_1\ell(|I|)\Big)\\
&\le
\sum_{I \in\mathcal{T}'} \mathbb{P}_{\theta}\Big(\sum_{i\in I_*\backslash  I} \xi_i^2 
\ge  M_\xi\ell(|I_*|)+(C_{A_1}-M_\xi)\ell(|I_*|)\Big)\\
&= \sum_{I\in[n]} \mathbb{P}_{\theta}\Big(\sum_{i\in I} \xi_i^2 
\ge  M_\xi \ell(|I|) + (C_{A_1}-M_\xi)\ell(|I_*|)\Big)  \notag\\
&\le H_0 e^{-\alpha'_1 \ell(|I_*|)},
\end{align*}
with $\alpha'_1=\alpha_\xi(C_{A_1}-M_\xi)$, $A_1$ is assumed to be so large that
$C_{A_1}> M_\xi$.
\end{proof}

\begin{proof} [Proof of  Theorem \ref{theorem_2}] 
First we prove \eqref{th2_i}. Consider the case $|I^*|\ge 1$. Let $B=\{|\tilde{I}|>M_1|I^*|\}$. 
Recalling that $I^*=I_*(A_0)$, by \eqref{property_Io}, $\theta_i^2 \le 
A_0 \sigma^2 \log(\tfrac{qn}{|I^*| \vee 1})$ for all $i\in (I^*)^c$. 
Using this, Condition \eqref{cond_nonnormal}, the definition \eqref{threshold}, 
property \eqref{th1_i'} of Theorem \ref{theorem_1}, we have  that for $I^*=I_*(A_0)$ 
and $\hat{I}=\hat{I}(K)$
\begin{align}
\label{relation_1} 
\textnormal{FPR}(\hat{I},I^*)&=\operatorname{E}_\theta\tfrac{|\hat{I}\backslash I^* |}{n-|I^*|}
=\tfrac{1}{n-|I^*|}\sum_{i\in (I^*)^c}\operatorname{E}_\theta \hat{\eta}(X_i)(\mathrm{1}_{B}
+\mathrm{1}_{B^c})
\notag\\
&\le  \mathbb{P}_\theta(B)+\tfrac{1}{n-|I^*|}
\sum_{i\in (I^*)^c}\mathbb{P}_\theta \big(\tfrac{2\theta_i^2}{\sigma^2}
+2\xi_i^2\ge K\log (\tfrac{qn}{M_1|I^*|})\big )\notag\\
&\le  \mathbb{P}_\theta(B)+\tfrac{1}{n-|I^*|}\sum_{i\in (I^*)^c}
\mathbb{P}_\theta \big(\xi_i^2\ge(\tfrac{K}{2}-A_0) \log(\tfrac{qn}{|I^*|})
-\tfrac{K}{2}\log M_1 \big)\notag\\
&\le 
H_0 e^{-\alpha_0\ell(|I^*|)} + H'(\tfrac{n}{|I^*|})^{-\alpha'},
\end{align} 
for $K$ and $A_0$ such that $\tfrac{K}{2}-A_0>C_\xi$ 
and the property \eqref{th1_i'} of Theorem \ref{theorem_1} can be applied, 
where  $\alpha'=\alpha_\xi(\tfrac{K}{2}-A_0)$ and  $H'=H_\xi e^{\alpha_\xi K \log \sqrt{M_1}}$.

Now we consider the case $|I^*|=0$ (i.e., $I^*=\varnothing$).
Reasoning similarly to \eqref{relation_1}, now with $B=\{|\tilde{I}|>\log n \}$, we derive
\begin{align}
\textnormal{FPR}(\hat{I},I^*) 
&\le \mathbb{P}_\theta( B)+\tfrac{1}{n}
\sum_{i\in [n]}\mathbb{P}_\theta \big(\tfrac{2\theta_i^2}{\sigma^2}
+2\xi_i^2\ge K\log (\tfrac{qn}{\log n})\big ) \notag\\
&\le  \mathbb{P}_\theta(B)+\tfrac{1}{n}\sum_{i\in [n]}
\mathbb{P}_\theta \big(\xi_i^2\ge(\tfrac{K}{2}-A_0) \log(qn)-\tfrac{K}{2}\log\log n \big) \notag\\
\label{relation_2}
&\le H'' n^{-\alpha''},
\end{align} 
The relations  \eqref{relation_1} and \eqref{relation_2}
establish \eqref{th2_i}.
 
Next, we proof the assertion \eqref{th2_ii}. If $I_*=\varnothing$, the claim follows, assume 
$|I_*|\ge 1$. Denote $B_\delta = \{|\tilde{I}|\le \delta |I_*|\}$. Using Condition \eqref{cond_nonnormal}, 
the definition \eqref{threshold}, \eqref{property_Io}, property \eqref{th1_ii'} 
of Theorem \ref{theorem_1}, and the fact that $(a+b)^2\ge 2a^2/3-2b^2$ for any 
$a,b\in\mathbb{R}$, we have  that for $I_*=I_*(A_1)$ and $\hat{I}=\hat{I}(K)$
\begin{align}
\label{NDR}
\textnormal{NDR}(\hat{I},I_*)&=\frac{1}{|I_*|}\operatorname{E}_\theta |I_*\backslash \hat{I}|
=\frac{1}{|I_*|}\sum_{i\in I_*}\operatorname{E}_\theta (1\!-\!\hat{\eta}(X_i))\!\notag\\
&=\frac{1}{|I_*|}\sum_{i\in I_*}\!\mathbb{P}_\theta \Big(|\theta_i+\sigma \xi_i|
<\sigma\big[K\log (\tfrac{qn}{|\tilde{I}|})\big]^{1/2}\Big) \notag\\
&\le\frac{1}{|I_*|}\sum_{i\in I_*}\mathbb{P}_\theta\Big(\tfrac{2}{3}\theta_i^2\!-\!2\sigma^2\xi_i^2
<\sigma^2K \log (\tfrac{qn}{|\tilde{I}|})\Big)\notag\\
&\le\frac{1}{|I_*|}\sum_{i\in I_*}
\mathbb{P}_\theta\Big(\xi_i^2>\tfrac{A_1}{3}\log (\tfrac{en}{|I_*|})
-\tfrac{K}{2}\log (\tfrac{qn}{|\tilde{I}|})\Big) \notag\\
&\le \mathbb{P}_\theta(B_\delta)+\frac{1}{|I_*|}\sum_{i\in I_*} 
\mathbb{P}_\theta\Big(\xi_i^2>\tfrac{A_1}{6}\log (\tfrac{qn}{|I_*|})
-\tfrac{K}{2}\log (\tfrac{qn}{|\tilde{I}|}), B_\delta^c \Big) \notag\\
&\le H_0  e^{-\alpha'_1 \ell(|I_*|)} +\frac{1}{|I_*|}\sum_{i\in I_*}
\mathbb{P}_\theta\big(\xi_i^2>(\tfrac{A_1}{6}-\tfrac{K}{2})\log (\tfrac{qn}{|I_*|})
+\tfrac{K}{2}\log \delta\big)\notag\\
&\le H_0 e^{-\alpha'_1 \ell(|I_*|)} +C_2(\tfrac{qn}{I_*})^{-\alpha'_3} 
\le H_2 (\tfrac{n}{|I_*|})^{-\alpha_3},
\end{align}
for sufficiently large $A_1$ (such that $\tfrac{A_1}{6}-\tfrac{K}{2}>M_\xi$ 
and the property \eqref{th1_ii'} of Theorem \ref{theorem_1} can be applied), 
where $\alpha'_3=\alpha_\xi(\tfrac{A_1}{3}-\tfrac{K}{2}-M_\xi)$ 
and  $C_2=H_\xi e^{-\alpha_\xi K \log \sqrt{\delta}}$. The relation 
\eqref{th2_ii} is proved.

Since $I^*=I_*(A_0(K))=I_*(A_1(K))=I_*$ for $\theta\in\Theta(K)$, the assertion 
\eqref{th2_iii} follows from the relations  \eqref{relation_1} and \eqref{NDR}:
uniformly in $\theta\in \Theta(K)$,
\begin{align*}
\operatorname{E}_\theta \big(|\hat{I}\backslash I_*|+|I_*\backslash \hat{I}|\big)
&\le (n-|I_*|)H_1 (\tfrac{n}{|I_*|\vee 1})^{-\alpha_2}
+ |I_*| H_2 (\tfrac{n}{|I_*|})^{-\alpha_3}\\
&\le H_3 n (\tfrac{n}{|I_*|\vee 1})^{-\alpha_4}.
\qedhere
\end{align*}
\end{proof}

\begin{proof}[Proof of Theorem \ref{theorem_3}]
First, we proof assertion \eqref{th3_i}. Introduce the event
$B_{\delta}=\{|\tilde{I}|< \delta |I_*|\}$. We argue along the same lines 
as in \eqref{relation_1} and \eqref{relation_2} for the two cases $|I_*|>0$ and $|I_*|=0$. 
For the case $|I_*|>0$, we use \eqref{th1_ii'} of Theorem \ref{theorem_1}, 
\eqref{th2_i} of Theorem \ref{theorem_2} and the fact that 
$I^*=I_*$ for $\theta\in \Theta(K)$ to derive 
\begin{align*}
\textnormal{FDR}(\hat{I}) &=\operatorname{E}_\theta 
\tfrac{|\hat{I}\backslash I_*|}{|\hat{I}|}
=\operatorname{E}_\theta\big[ \tfrac{|\hat{I}\backslash I^*|}{|\hat{I}|}
(\mathrm{1}_{B_{\delta}^c}+\mathrm{1}_{B_{\delta}})\big]
\le \tfrac{1}{\delta|I_*|}\operatorname{E}_\theta |\hat{I}\backslash I^*|+\mathbb{P}_\theta (B_\delta) \\
&\le \tfrac{1}{\delta|I_*|} H_1n \big(\tfrac{n}{|I^*|\vee 1}\big)^{-\alpha_2} +
H_0 e^{-\alpha'_1\ell(|I_*|)}
\le
H_5 (\tfrac{n}{|I_*|})^{-\alpha_6}.
\end{align*}
The case $|I_*|=0$ is handled as follows. If $|\hat{I}|=0$, the claim holds. Assume 
$|\hat{I}|\ge 1$, then 
\begin{align*}
\textnormal{FDR}(\hat{I}) 
\le \operatorname{E}_\theta |\hat{I}\backslash I^*| 
\le H_1 n n^{-\alpha_2} = H_1 n^{-(\alpha_2-1)}.
\end{align*}
	
Next, we prove assertion \eqref{th3_ii}.  Introduce the event $B=\{|\tilde{I}|>M_1 |I^*|\}$.
Consider two cases: the case $|I^*| \ge n/(2M_1)$ and the case $|I^*|<n/(2M_1)$. 
Suppose $|I^*| \ge n/(2M_1)$, then 
$\textnormal{FNR}(\hat{I})=\operatorname{E}_\theta \tfrac{|I_*\backslash \hat{I}|}{n-|\hat{I}|}
\le 1 \le 2 M_1 \tfrac{|I^*|}{n}=2 M_1 \tfrac{|I_*|}{n}$ and \eqref{th3_ii} holds, as 
$I^*=I_*$ for $\theta\in\Theta(K)$.

Now suppose  $|I^*|<n/(2M_1)$. Then using the same reasoning as in \eqref{NDR}, 
Theorem \ref{theorem_1} and the fact that $I^*=I_*$ for $\theta\in\Theta(K)$,
we again obtain \eqref{th3_ii}:
\begin{align*}
\textnormal{FNR}(\hat{I})&=\operatorname{E}_\theta \tfrac{|I_*\backslash \hat{I}|}{n-|\hat{I}|}
=\operatorname{E}_\theta\big[ \tfrac{|I_*\backslash \hat{I}|}{(n-|\hat{I}|)}
(\mathrm{1}_{B^c}+\mathrm{1}_{B})\big]\\&\le\tfrac{1}{n-M_1|I^*|} 
\operatorname{E}_\theta |I_*\backslash \hat{I}|+ \mathbb{P}_\theta(B)
\le H_6 (\tfrac{n}{|I_*|})^{-\alpha_7}. \hfil  \qedhere
\end{align*}	
\end{proof}

\begin{proof}[Proof of  Theorem \ref{theorem_uncertainty}]
We first establish the coverage property. 
Recall that $\hat{r}=\hat{r}(\tilde{I})= n (\tfrac{n}{|\tilde{I}|\vee 1})^{-\alpha'_4}$ and by 
\eqref{th2_iii} of Theorem \ref{theorem_2},
\begin{align}
\label{th4_rel1}
\sup_{\theta\in\Theta(K)} \operatorname{E}_\theta | \hat{\eta} - \eta_*| 
&\le H_3 n (\tfrac{n}{|I_*|\vee 1})^{-\alpha_4}. 
\end{align}
Denote $B_{\delta}=\{|\tilde{I}|\le \delta |I_*|\}$.
Now, by using the Markov inequality, property \eqref{th1_ii'} of Theorem \ref{theorem_1} 
and \eqref{th4_rel1}, we obtain, uniformly in $\theta\in\Theta(K)$,
\begin{align*}
\mathbb{P}_{\theta}\big(\eta_*\notin B(\hat{\eta},\hat{r})\big) 
&= \mathbb{P}_{\theta}\big(| \eta_*-\hat{\eta}|> \hat{r}, B^c_{\delta} \big)+
\mathbb{P}_{\theta}(B_{\delta}) \\
&\le 
\mathbb{P}_{\theta}(B_{\delta})+\mathbb{P}_{\theta}\big(| \eta_*-\hat{\eta}|
> n (\tfrac{n}{ (\delta |I_*|) \vee 1})^{-\alpha'_4} \big)\\
&\le  H_0 e^{-\alpha'_1 \ell(|I_*|)} + H_3 \delta^{\alpha'_4} 
(\tfrac{n}{|I_*| \vee 1})^{-(\alpha_4- \alpha'_4)}
\le H_7 (\tfrac{n}{|I_*| \vee 1})^{-\alpha_7},
\end{align*}
which proves the  coverage property.

It remains to prove the size property. Let $B= \{|\tilde{I}|\ge M_1|I^*|\}$.
For any $M'_1> M_1^{\alpha'_4}\vee 1$, we have that, 
uniformly in $\theta\in\Theta(K)$ (so that $I^*=I_*$),
\begin{align*}
\mathbb{P}_{\theta}\big(\hat{r}\ge M'_1 r_*, B^c\big)
\le \mathbb{P}_{\theta} 
\Big(n \big(\tfrac{n}{M_1 |I^*|\vee 1}\big)^{-\alpha'_4} 
\ge M'_1 n \big(\tfrac{n}{|I_*|\vee 1}\big)^{-\alpha'_4}\Big) =0.
\end{align*}
Using this, the property \eqref{th1_i'} of Theorem \ref{theorem_1} and the fact that 
$I^*=I_*$ for $\theta\in\theta(K)$, we derive that, uniformly in $\theta\in\Theta(K)$,
\begin{align*}
\mathbb{P}_{\theta}\big(\hat{r}\ge M'_1 r_*\big)
&\le \mathbb{P}_{\theta}(B)+
\mathbb{P}_{\theta}\big(\hat{r}\ge  M'_1 r_*, B^c\big) \\
&=  \mathbb{P}_{\theta}(B) \le H_0 e^{-\alpha_0\ell(|I^*|)}  
\le  H_8 (\tfrac{n}{|I^*| \vee 1 })^{-\alpha_8},
\end{align*}
yielding the size property.
\end{proof}

\end{document}